\documentclass[11pt]{article}

\usepackage{tgtermes}
\usepackage[T1]{fontenc}
\usepackage{amsmath,amsthm}
\usepackage{graphicx}
\usepackage{xcolor}

\usepackage{arxivs}
\usepackage[utf8]{inputenc}
\usepackage{amsfonts,amssymb}

\usepackage{cite}
\usepackage{hyperref}
\usepackage{textcomp}
\usepackage{tabularx}
\usepackage{caption}
\usepackage{subcaption}
\usepackage{diagbox}

\newcommand{\Eb}{\mathbb{E}}
\newcommand{\Rb}{\mathbb{R}}
\newcommand{\bR}{\bar{\mathbb{R}}}  
\newcommand{\one}{\mathbb{I}}
\newcommand{\onen}{\mathbf{1}}

\newcommand{\Ac}{\mathcal{A}}
\newcommand{\At}{\tilde{\mathcal{A}}}
\newcommand{\Fc}{\mathcal{F}}
\newcommand{\Kc}{\mathcal{K}}
\newcommand{\Lc}{\mathcal{L}}
\newcommand{\Sc}{\mathcal{S}}
\newcommand{\Zc}{\mathcal{Z}}
\newcommand{\Zp}{\mathcal{Z}_p}
\newcommand{\Lpm}{\mathcal{L}_p(\Omega, \Fc, P; \Rb^m)}  
\newcommand{\Lqm}{\mathcal{L}_q(\Omega, \Fc, P; \Rb^m)}
\newcommand{\X}{\mathcal{X}}
\newcommand{\Y}{\mathcal{Y}}
\newcommand{\D}{\mathcal{D}}

\DeclareMathOperator{\dom}{dom}
\DeclareMathOperator{\var}{VaR}
\DeclareMathOperator{\avar}{AVaR}
\DeclareMathOperator{\Mvar}{MVaR}   
\DeclareMathOperator{\Mavar}{MAVaR}  
\DeclareMathOperator{\Vmavar}{VMAVaR} 
\DeclareMathOperator{\argmax}{\arg\max}

\newtheorem{definition}{Definition}
\newtheorem{theorem}{Theorem}
\newtheorem{corollary}{Corollary}
\newtheorem{proposition}{Proposition}

\title{On Risk Evaluation and Control of Distributed Multi-Agent Systems}
\author{Aray Almen \\
		Department of Mathematical sciences \\
		Stevens Institute of Technology \\
		Hoboken, NJ 07030, USA \\
		\texttt{aalmen@stevens.edu} \\
		\And 
		Darinka Dentcheva \\
		Department of Mathematical sciences \\
		Stevens Institute of Technology \\
		Hoboken, NJ 07030, USA \\
		\texttt{darinka.dentcheva@stevens.edu}
		}

\date{}

\begin{document}
\maketitle

\begin{abstract}
	In this paper, we deal with risk evaluation and risk-averse optimization of complex distributed systems with general risk functionals. We postulate a novel set of axioms for the functionals evaluating the total risk of the system.  We derive a dual representation for the systemic risk measures and propose a way to construct non-trivial families of measures by using either a collection of linear scalarizations or non-linear risk aggregation. The new framework facilitates risk-averse sequential decision-making by distributed methods. The proposed approach is compared theoretically and numerically to some of the systemic risk measurements in the existing literature. 

	We formulate a two-stage decision problem with monotropic structure and systemic measure of risk. The structure is typical for distributed systems arising in energy networks, robotics, and other practical situations. A distributed decomposition method for solving the two-stage problem is proposed and it is applied to a problem
	arising in communication networks. 
	We have used this problem to compare the methods of systemic risk evaluation. We show that the proposed risk aggregation leads to less conservative risk evaluation and results in a substantially better solution of the problem at hand as compared to an aggregation of the risk of individual agents and other methods.
 
\end{abstract}


\section{Introduction}
\label{sec:introduction}

	Evaluation of the risk of a system consisting of multiple agents is one of the fundamental problems relevant to many fields. A crucial question is the assessment of the total risk of the system taking into account the risk of each agent and its contribution to the total risk. Another issue arises when the risk evaluation is based on confidential or proprietary information. There is extensive literature addressing the properties of risk measures and their use in finance. Our goal is to address situations related to robotics, energy systems, business systems, logistic problems, etc. The analysis in financial literature may not be applicable in such situations due to the heterogeneity of the sources of risk, the nature, and the complexity of relations in those systems.  In many systems, the source of risk is associated with highly non-trivial aggregation of the features of its agents, which may not be available in an analytical form. For example, in automated robotic systems, the exchange of information may be limited or distorted due to the speed of operation, the distance in space between the agents, or other reasons. Another difficulty associated with the evaluation of risk arises when the risk of one agent stems from various sources of uncertainty of different nature. The question of how to aggregate those risk factors in one loss function does not have a straightforward answer. 

	The risk of one loss function can be evaluated using a coherent measure of risk such as Average Value-at-Risk, mean-semideviation or others. More traditional (non-coherent) measures of risk such as Value-at-Risk (VaR) are also very popular and frequently used. We refer to \cite{Follmer} for an extensive treatment of risk measures for scalar-valued random variables, as well as to \cite{textbook} where risk-averse optimization problems are analyzed as well.

	The main objective of this paper is to suggest a new approach to the risk of a distributed system and show its viability and potential in application to risk-averse decision problems for distributed multi-agent systems. While building on the developments thus far, our goal is to identify a framework that is theoretically sound but also amenable to efficient numerical computations for risk-averse optimization of large multi-agent systems.
	We propose a set of axioms for functionals defined on the space of random vectors. The random vector is comprising risk factors of various sources, or is representing the loss of each individual agent in a multi-agent system. While axioms for random vectors have been proposed earlier, our set of axioms differs from those in the literature most notably with respect to the translation equivariance condition, which we explain in due course. The resulting systemic risk measures reduce to coherent measures of risk for scalar-valued random variables when the dimension of the random vectors becomes one. We derive the dual representation of the systemic measures of risk with less assumptions than known for multi-variate risks. In our derivation, we establish one-to-one correspondences between the axioms and properties of the dual variables. We also propose several ways to construct systemic risk measures  and analyze their properties. The important features of the proposed measures are the following. They are conformant with the axioms; they can be calculated efficiently, and are amenable to distributed optimization methods. 

	We have formulated a risk-averse two-stage optimization problem with a structure, which is typical for a system of loosely coupled subsystems. The proposed numerical method is applied to manage the risk of a distributed operation of agents. The distributed method lets each subsystem optimize its operation with minimal information exchange among the other subsystems (agents). This aspect is important for multi-agent systems where some proprietary information is involved or when privacy concerns exist. The method demonstrates that distributed calculation of the systemic risk is possible without a big computational burden. We then consider a two-stage model in wireless communication networks, which extends the static model discussed in \cite{Ma-me-Zavlanos}. It addresses a situation when a team of robots explores an area and each robot reports relevant information. The goal is to determine a few reporting points so that the communication is conducted most efficiently while managing the risk of losing information.   We conduct several numerical experiments to compare various systemic risk measures. 

	Our paper is organized as follows. In section \ref{s:prelim} we provide preliminary information on coherent measures of risk for scalar-valued random variables and survey existing methods for risk evaluation of complex systems. Section \ref{s:axiomatic} contains the set of axioms, the dual representation associated with the resulting systemic risk measures, and two ways to construct such measures in practice. Section \ref{s:comparison-theoretical} provides a theoretical comparison of the new measures of risk to other notions. In particular, we discuss other sets of axioms, explore relations to two notions of multivariate Average Value-at-Risk, and pay attention to the effect of the aggregation of risk before and after risk evaluation. 
	In section \ref{s:numerical}, we formulate a risk-averse two-stage stochastic programming problem modeling wireless information exchange and seeking to locate a constraint number of information exchange points. We devise a distributed method for solving the problem and report a numerical comparison with several measures of risk, and other systemic measures. We pay attention to the comparison between the principles of aggregation for the purpose of total risk evaluation.

\section{Preliminaries}
\label{s:prelim}

% --------------------------------- COHERENT RISK MEASURES ------------------------------ %
\subsection{Coherent risk measures}
\label{s:prelim-scalar-measures}

	The widely accepted axiomatic framework for coherent measures of risk was proposed in \cite{Artzner} and further analyzed in \cite{Delbaen}, \cite{Follmer}, \cite{Leitner:2005}, \cite{RS:2005,RuSh:2006a}, \cite{PflRom:07} and many others works.
	It is worth noting that another axiomatic approach was initiated in \cite{KijOhn:1993} and this line of thinking was developed to an entire framework in \cite{rockafellar2013fundamental}.
	For a detailed exposition, we refer to \cite{textbook} and the references therein. 
	Let $\Lc_p(\Omega, \Fc, P)$ be the space of real-valued random variables, defined on the probability space $(\Omega,\Fc, P)$, that have finite $p$-th moments, $p \in [1, \infty)$, and are indistinguishable on events with zero probability.  We shall assume that the random variables represent random costs or losses. A lower semi-continuous functional $\varrho: \Lc_p(\Omega, \Fc, P) \to \Rb\cup\{+\infty\}$ is a \textit{coherent risk measure} if it is convex, positively homogeneous, monotonic with respect to the a.s. comparison of random variables, and satisfies the the following translation property
\[
\varrho[X + a] = \varrho[X] + a \text{ for all } X \in \Lc_p(\Omega, \Fc, P),\; a \in \Rb. 
\]
	If $\varrho[\cdot]$ is monotonicity, convex, and satisfies the translation property, then it is called a \emph{convex} risk measure.  Some examples of coherent measures of risk include Average Value-at-Risk (also called Conditional Value-at-Risk) and mean-semideviations measure, which are defined as follows. The Average Value-at-Risk at level $\alpha$ for a random variable $Z$ is defined as 
	\[
	\avar_\alpha [Z] = \inf_{\eta \in \Rb} \Big\{ \eta + \frac{1}{\alpha} \Eb \big[ (Z - \eta)_+ \big] \Big\}
	\]
	It is a special case of the higher-order measures of risk:
	\[
	\varrho[Z] = \min_{t\in \Rb}\bigg\{ t + \frac{1}{\alpha} \big\| \big(Z-t\big)_+\big\|_p \bigg\},\quad \alpha \in (0,1),
	\]
	where $\|\cdot\|_p$ refers to the norm in $\Lc_p(\Omega, \Fc, P)$. 
	The mean semi-deviation of order $p$ is given by
	\[
	\varrho[Z] = \Eb[Z] + \varkappa \big\| \big(Z-\Eb[Z]\big)_+\big\|_p,\quad \varkappa \in [0,1].
	\] 
	The space $\Lc_p(\Omega, \Fc, P)$ equipped with its norm topology is paired with the space $\Lc_q(\Omega, \Fc, P)$ equipped with the weak$^*$ topology where $\frac{1}{p} + \frac{1}{q} = 1$. For any $Z \in \Zc $ and $\xi \in \Zc^*$, we use the bilinear form:
	\[ 	
	\langle \xi, Z \rangle = \int_\Omega \xi(\omega) Z(\omega) dP(\omega). 
	\]
	The following result is known as a dual representation of coherent measures of risk. A proper lower semicontinuous coherent risk measure $\varrho$ has a dual representation
	\begin{equation}
    	\varrho[Z] = \sup_{\xi \in \Ac_\varrho} \langle \xi, Z \rangle, \quad Z \in \Zc,
    	\label{dual rep rho}
	\end{equation}
	where $\Ac_\varrho \subset \{ \xi \in \Zc^* ~ | ~ \xi \geq 0, ~ \int_\Omega \xi(\omega)P(d\omega) = 1 \}$ is the convex-analysis subdifferential  $\partial \varrho[0]$.

	Risk measures have also been defined by specifying a set of desired values for the random quantity in question; this set is called an acceptance set. Denoting the acceptance set by $\Kc\subset\Rb$, the risk of a random outcome $Z$ is defined as:
	\begin{equation}
    	\varrho_{\Kc}[X] = \inf \{ m \in \Rb ~|~ X - m \in \Kc \}.
    	\label{acceptance set rho}
	\end{equation}
	In finance, this notion of risk is interpreted as the minimum amount of capital that needs to be invested to make the final position acceptable. It is easy to verify that $\varrho [\cdot]$ in \eqref{acceptance set rho} is a coherent measure if and only if $\Kc$ is a convex cone (cf. \cite{FollmerSchied}).

% ------------------------- RISK MEASURES FOR COMPLEX SYSTEMS ------------------------- %

\subsection{Risk measures for complex systems}
\label{s:prelim-review}

	As the risk is not additive, when we deal with distributed complex systems, we need to address the question of risk evaluation for the entire system. This risk is usually called systemic in financial literature and the proposed measures for its evaluation are termed systemic risk measures. 

	Assume that the system consists of $m$ agents. One approach to evaluating the risk of a system is to use an \textit{aggregation function}, $\Lambda:\Rb^m\to \Rb$, and univariate risk measures. Let $X \in \Lpm$ be an $m$-dimensional random vector comprising the costs incurred by the system, where each component $X_i$ corresponds to the costs of one agent. The first approach to systemic risk is to choose a univariate risk measure $\varrho_0$ and apply it to the aggregated cost $\Lambda(X).$ If we prefer to use an acceptance set $\Kc$ as in \eqref{acceptance set rho}, the systemic risk can be defined as:
	\begin{align}
	    \varrho[X] = \varrho_0[\Lambda(X)] = \inf \big\{ z \in \Rb ~|~ \Lambda(X) - z \in \Kc \big\}. 
	    \label{def: risk of lambda}
	\end{align}
	In (\cite{Chen} this point of view is analyzed in finite probability spaces and it is shown that any monotonic, convex, positively homogeneous function provides a risk evaluation as in \eqref{def: risk of lambda} as long as it is consistent with the preferences represented in the definition of $\Kc$. The point of view presented in definition \eqref{def: risk of lambda}  is further extended in \cite{Kromer}, where the authors analyzed convex risk measures defined on a general measurable space and proposed examples of aggregation functions suitable for a financial system. In both studies, the structural decomposition of the systemic risk measure \eqref{def: risk of lambda} is established when the aggregation function $\Lambda$ satisfies properties similar to the axioms postulated for risk measures. In \cite{Brunnermeier}, the authors considered a particular case of an aggregation function, proposing an evaluation method for the risk associated with the cumulative externalities or costs endured by financial institutions. Note that these evaluation methods rely on a choice of one aggregation function suitable for a specific problem. 
	\\
	The translation property for constant vectors is introduced in \cite{Burgert} for convex risk measures defined for bounded random vectors. This property differs from the one we propose here. The authors analyzed the maximal risk over a class of aggregation functions rather than using one specific function. We refer to \cite{RuschendorfBook} for an overview of the risk measures constructed this way. A similar approach is taken in \cite{ekeland2011law}, where law-invariant risk measures for bounded random vectors are investigated for the purpose of obtaining a Kusuoka representation. The axioms proposed in \cite{Burgert,ekeland2011law} are closest to ours and we provide more detailed discussion in section \ref{s:axiomatic}.

	Another approach to risk evaluation of complex systems consists of evaluation of the risk of individual agents first and aggregation of the obtained values next. This method is used, for example, in \cite{Biagini} and \cite{Feinstein}. Using the notion of acceptance sets the systemic risk measure is defined in \cite{Biagini} in the following way:
	\begin{align*}
	    \varrho[X] & = \varrho_0[\Lambda(X)] = \inf \Big\{ \sum_{i=1}^m z_i ~|~ z \in \Rb^m, ~\Lambda(X - z) \in\Kc\subset\Rb^m \Big\}. 
	\end{align*}
	The proposed measures of risk in section \ref{s:axiomatic} also accommodate this point of view. 
	A further extension in \cite{Biagini} replaces the constant vector $z \in \Rb^m$ by a random vector $Y \in \mathcal{C}$, where $\mathcal{C}$ is a given set of admissible allocations. This formulation of the risk measure allows to decide scenario-dependent allocations, where the total amount $\sum_{i=1}^m z_i$ can be determined ahead of time while individual allocations $z_i$ may be decided in the future when uncertainty is revealed.  In \cite{Feinstein} a set-valued counterpart of this approach is proposed by defining the systemic risk measure as the set of all vectors that make the outcome acceptable. Once the set of all acceptable allocations is constructed, one can derive a scalar-valued \textit{efficient allocation rule} by minimizing the weighted sum of components of the vectors in the set. 
	Set-valued risk measures were proposed in \cite{jouini2004vector}, see also \cite{Ararat,hamel2010duality} for duality theory including the dual representation for certain set-valued risk measures. 
	In fast majority of literature, the systemic risk depends on the choice of the aggregation function $\Lambda$ and how well it captures the interdependence between the components. 
	To capture the dependence, an approach based on copula theory was put forward in \cite{pflug2018systemic}. It is assumed that independent operation does not carry systemic risk and, hence, the local risk can be optimized by each agent independently. The systemic risk measures are then constructed based on the copulas of the distributions.

	Another line of work includes methods that use some multivariate counterpart of the univariate risk measures. The main notion here is the Multivariate Value-at-Risk ($\Mvar$) for random vectors, which is identified with the set of \textit{$p$-efficient points}. Let $F_X(\cdot)$ be the right-continuous distribution function of a random vector $X$ with realizations in $\Rb^m$. A $p$-efficient points for $X$ is a point $v\in\Rb^m$ such that $F_X(v)\geq p$ and there is no point $z$ that satisfies $F_X(z)\geq p$ with $z\leq v$ componentwise. This notion plays a key role in optimization problems with chance constraints (see e.g. \cite{textbook}). Multivariate Value-at-Risk satisfies the properties of translation equivariance, positive homogeneity and monotonicity.  This notion is used to define Average Value-at-Risk for multivariate distributions ($\Mavar$) in \cite{LeePrekopa,Noyan2013optimization,Prekopa}.
	Let ${Z}_p$ be the set of all points, each of which is component-wise larger than some $p$-efficient point:
	\[
	Z_p = \bigcup_{s \in \Mvar_p(X)} (s + \Rb^m_+). 
	\]
	In \cite{LeePrekopa}, Lee and Prekopa define the $\Mavar$ of a random vector $X$ at level $p$ as
	\begin{equation}
	    \Mavar_p(X) = \Eb (\Lambda(X) ~|~ X \in D_p \}, \label{Prek_Mavar}
	\end{equation}
	where $\Lambda$ is assumed integrable with respect to $F_X$, i.e., $\Eb(\Lambda(X))$ is finite.  
	It is shown in \cite{LeePrekopa} that $\Mavar$ is translation equivariant, positive homogeneous and subadditive only when all of the components of the random vector are independent. 

	While the definition of $\Mavar$ above is scalar-valued, in \cite{Merakli} the authors define a Multivariate Average Value-at-Risk ($\Mavar$) using the notion of $p$-efficient points as $\Mvar_p(X)$ and the extremal representation of the Average Value-at-Risk. First for given probability $p\in(0,1)$, we consider the vectors
	\[
	    \Mavar_p(X; v) = v + \frac{1}{p} \Eb[(X-v)_+], 
	\] 
	where $ [(X-v)_+]_i = \max(0, X_i - v_i) $, $i = 1,\dots,m$.
	Then, the following vector-optimization problem is solved:
	\[
	    \Vmavar_p(X) = \min \{ \Mavar_p(X; v) : v\in \Mvar_p(X) \}.
	\]
	The vector-valued Multivariate Average Value-at-Risk is monotonic, positively homogeneous, translation equivariant, but is not subadditive. Note that in both $\Mvar$ and $\Mavar$, one needs to use a scalarization function to obtain a scalar value for the risk.

	We shall compare our proposal to the aforementioned risk measures in section \ref{s:comparison-theoretical}.

% --------------------------------- AXIOMATIC APPROACH ---------------------------------- %

\section{Axiomatic Approach to Risk Measures for Random Vectors}
\label{s:axiomatic}

	In this section, we propose a set of axioms to measures of risk for random vectors with realizations in $\Rb^m$. This framework is analogous to the coherent risk measures properties for scalar-valued random variables. In fact, if $m=1$, the proposed set of axioms exactly coincides with those in \cite{textbook}. 
    We denote by $\Zc=\Lc_p(\Omega, \Fc, P;\Rb^m)$ be the space of random vectors with realizations in $\Rb^m$, defined on $(\Omega,\Fc, P)$. 
    Throughout the paper, we shall consider risk measure $\varrho$ for random vectors in $\Zc$ to be a lower-semi-continuous functional $\varrho: \Zc \to \Rb\cup\{+\infty\}$ with non-empty domain. We denote the $m$-dimensional vector, whose components are all equal to one by $\onen$ and the 
	random vector with realizations equal to $\onen$ by $\one$.
	\begin{definition}
	\label{d:riskonvectors}
	    A lower semi-continuous functional $\varrho: \Zc \to \Rb\cup\{+\infty\}$ is a coherent risk measure with preference to small outcomes, iff it satisfies the following axioms:
	    \begin{itemize}
	        \item[A1.] Convexity: For all $X, Y \in \Zc$ and $ \alpha \in (0,1)$, we have:
	            $$ ~ \varrho[\alpha X + (1-\alpha)Y] \leq \alpha \varrho[X] + (1 - \alpha)\varrho[Y]. $$
	        \item[A2.] Monotonicity: For all $X, Y \in \Zc$, if $X_i \geq Y_i$ for all components $i = 1, \dots, m$ $P$-a.s., then $\varrho[X] \geq \varrho[Y]$.
	        \item[A3.] Positive homogeneity: For all $X \in \Zc$ and $t > 0$, we have $\varrho[tX] = t\varrho[X]$.
	        \item[A4.] Translation equivariance: For all $X \in \Zc$ and $a \in \Rb$, we have $\varrho[X + a\one] = \varrho[X] + a\varrho[\one].$
	    \end{itemize}
	    A lower semi-continuous  functional $\varrho: \Zc \to \Rb\cup\{+\infty\}$ is a convex risk measure with preference to small outcomes, iff it satisfies axioms A1, A2, and A4.
	\end{definition}

	The axioms of convexity and positive homogeneity are defined in a similar way to the properties of coherent risk measures, while the random vectors are now compared component-wise for the property of monotonicity. 
	% The monotonicity property implies that if all components of a random vector incur a smaller loss a.s., then the risk of this outcome should be smaller as well. The convexity property ensures that the diversification of random loss leads to smaller values of risk. The positive homogeneity implies that if the random cost increases by a constant factor of $t>0$, then the risk should also increase proportionally. 
	The main difference is the definition of a translation equivariance axiom. It suggests that if the random loss increases by a constant amount for all components, then the risk should also increase by the same amount. These axioms differ from the previous axioms proposed in the literature. 

	\subsection{Dual representation}

	In order to derive a dual representation of the multivariate risk measure, we pair the space of random vectors $\Zc = \Lpm$, $p \in [1, \infty)$ with the space $\Zc^* = \Lqm$, where $q \in (1, \infty)$ is such that $\frac{1}{p} + \frac{1}{q} = 1$, $q = \infty$ for $p = 1$. For $X \in \Zc$ and $\zeta \in \Zc^*$ the bilinear form $\langle \cdot, \cdot \rangle $ on the product space $\Zc \times \Zc^*$ is defined as follows:
	\[ 
	\langle \zeta, X \rangle = \int_\Omega \langle \zeta(\omega), X(\omega) \rangle dP(\omega). 
	\]
	The Fenchel conjugate function $\varrho^* : \Zc^* \to \Rb\cup\{+\infty\}$ of the risk measure $\varrho$ is given by
	\[
	\varrho^*[\zeta] = \sup_{X \in \Zc} \big\{ \langle \zeta, X \rangle - \varrho[X] \big\}, 
	\]
	and the conjugate of $\varrho^*$ (the bi-conjugate function) is
	\[
	\varrho^{**}[X] = \sup_{\zeta \in \Zc^*} \big\{ \langle \zeta, X \rangle - \varrho^*[\zeta] \big\}. 
	\]
	Fenchel-Moreau theorem implies that if $\varrho[\cdot]$ is convex and lower semicontinuous, then $\varrho^{**} = \varrho$ and that 
	\begin{equation}
	    \varrho[X] = \sup_{\zeta \in \At} \big\{ \langle \zeta, X \rangle - \varrho^*[\zeta] \big\}
	    \label{general dual}
	\end{equation}
	where $\At = \dom(\varrho^*)$ is the domain of the conjugate function $\varrho^*$. Then based on the Fenchel-Moreau theorem and the axioms proposed in this paper, we show the following theorem.
	\begin{theorem}
	\label{dual-rep}
	    Suppose $\varrho: \Zc \to \bR\cup\{+\infty\}$ is a convex and lower semicontinuous risk functional. Then the following holds:
	    \begin{itemize}
	        \item[(i)] Property A2 is satisfied if and only if $\zeta \geq 0$ a.s. for all $\zeta$ in the domain of $\varrho^*$.
	        \item[(ii)] Property A3 is satisfied if and only if $\varrho^*$ is the indicator function of $\partial\varrho[0]$, i.e.
	            \begin{equation}
	                \varrho[X] = \sup_{\zeta \in \Ac} \{ \langle \zeta, X \rangle \}.
	                \label{dual}
	            \end{equation} 
	            \item[(iii)] Property A4 is satisfied if and only if $\varrho[\one ] = \langle \one, \mu_\zeta \rangle $ for all $\zeta \in \Ac$, where $\mu_\zeta = \int_\Omega \zeta(\omega) P(d\omega)$. 
	    \end{itemize}
	\end{theorem}

	\begin{proof} 
	    Since $\varrho[\cdot]$ is convex and lower semicontinuous and we have assumed that it has non-empty domain, the representation \eqref{general dual} holds by virtue of the Fenchel-Moreau theorem. 

	    (i) Suppose $\varrho$ satisfies the monotonicity condition. Assume that $\zeta_i(\omega)< 0 $ for $\omega \in \Delta \in \Fc$ with $P(\Delta) > 0$ for some component $i = 1,\dots,m$. Define $\bar{X}_i$ equal to the indicator function of the event $\Delta$ and $\bar{X}_j = 0 $ for $j \neq i, j = 1,\dots, m$. Take any $X$ with support in $\Delta$ such that $\varrho[X]$ is finite and define $X_t: = X - t\bar{X}$. Then for $t \geq 0$, we have that $X \succeq X_t$, and $\varrho[X] \geq \varrho[X_t]$ by monotonicity. Consequently,
	    \begin{align*}
	        \varrho^*[\zeta] &\geq \sup_{t \in \Rb_+} \Big\{ \langle \zeta, X_t \rangle - \varrho[X_t] \Big\} \geq \sup_{t\in \Rb_+} \Big\{\langle \zeta, X \rangle - t\langle\zeta,\bar{X}\rangle -\varrho[X] \Big \} \\
	        &= \sup_{t\in \Rb_+} \Big\{\langle \zeta, X \rangle - t \int_\Delta \zeta_i(\omega)P(d\omega) -\varrho[X] \Big\} = +\infty.
	    \end{align*}
	    It follows that $\varrho^*[\zeta] = +\infty$ for every $\zeta \in \Zc^*$ with at least one negative component, thus $\zeta \notin \dom\varrho^*$. Conversely, suppose that $\zeta\in\Zc^*$ has realizations in $\Rb^m$ with nonnegative components $P$-a.s. Then whenever $X \succeq X'$, we have:
	    \begin{align*}
	        \langle \zeta, X \rangle & = \int_\Omega \langle \zeta (\omega), X(\omega) \rangle dP(\omega) \geq \int_\Omega \langle \zeta (\omega), X'(\omega)\rangle dP(\omega)=  \langle \zeta, X' \rangle.
	    \end{align*}
	    Consequently,
	    \begin{align*}
	        \varrho[X] & = \sup_{\zeta \in \Zc^*} \Big\{ \langle \zeta, X \rangle - \varrho^*[\zeta]\Big\} \geq \sup_{\zeta \in \Zc^*} \Big\{ \langle \zeta, X' \rangle - \varrho^*(\zeta) \Big\} = \varrho[X'].
	    \end{align*}
	    Hence, the monotonicity condition holds.

	    (ii) Suppose the positive homogeneity property holds. If $\varrho[tX] = t\varrho[X]$ for all $X \in \Zc$, then for any fixed $t>0$, we get
	    \begin{multline*}
	        \varrho^*[\zeta] = \sup_{X \in \Zc} \Big\{ \langle \zeta, X \rangle - \varrho[X] \Big\}
	        = \sup_{X \in \Zc} \Big\{ \langle \zeta, tX \rangle - \varrho[tX] \Big\} = 
	        \sup_{X \in \Zc} t \Big\{ \langle \zeta, X \rangle - \varrho[X] \Big\} =
	         t\varrho^*[\zeta].
	    \end{multline*}
	    Hence, if $\varrho^*[\zeta]$ is finite, then $\varrho^*[\zeta] = 0$ as claimed. 
	    Conversely, if $\varrho[X] = \sup_{\zeta \in \dom \varrho^*}\langle \zeta, X \rangle$, then $\varrho$ is positively homogeneous as a support function of a convex set.
	 
	    (iii) Suppose the translation property is satisfied, i.e. \\
	    $\varrho[X + t\one] = \varrho[X] + t\varrho[\one]$ for any $X \in \Zc$ and a constant $t \in \Rb$. Then for any $k \in \Rb$ and $\zeta \in \Zc^*$, we get:
	    \begin{align*}
	        \varrho^*[\zeta] &= \sup_{X \in \Zc} \Big\{ \langle \zeta, X + k \one \rangle - \varrho[X+k\one] \Big\} = \sup_{X \in \Zc} \Big\{ \int_\Omega \langle \zeta(\omega), X(\omega) + k\one \rangle P(d\omega) - \varrho[X] - k\varrho[\one] \Big\} \\
	        & = \sup_{X \in \Zc} \Big\{ \langle \zeta, X \rangle + k \int_\Omega \langle \one, \zeta(\omega) \rangle P(d\omega) - \varrho[X] - k\varrho[\one] \Big\} = \varrho^*[\zeta] + k \Big( \int_\Omega \langle \one, \zeta(\omega) \rangle P(d\omega) - \varrho[\one] \Big).
	    \end{align*}
	    
	    If $\varrho^*[\zeta]$ is finite, then $\varrho[\one] = \int_\Omega \langle \one, \zeta(\omega) \rangle P(d\omega)$. 

	    Let us denote $\mu_\zeta = \int_\Omega \zeta(\omega)P(d\omega)$, then we obtain  
	    \begin{equation}
	    \label{rhoofone}
	    \varrho[\one] = \langle \one, \int_\Omega \zeta(\omega) P(d\omega) \rangle = \langle \one, \mu_\zeta \rangle
	    \quad\text{for all }\zeta \in \Ac.
	    \end{equation}
	    Conversely, suppose $\varrho[\one] = \langle \one, \mu_\zeta \rangle$. Then for any $X \in \Zc$ and $k \in \Rb$:
	    \begin{align*}
	        \varrho[X + k\one] &= \sup_{\zeta \in \Zc^*} \Big\{ \langle \zeta, X + k \one \rangle - \varrho^*[\zeta] \Big\} = \sup_{\zeta \in \Zc^*} \Big\{ \int_\Omega \langle \zeta(\omega), X(\omega) + k\one \rangle P(d\omega) - \varrho^*[\zeta] \Big\} \\
	        & = \sup_{\zeta \in \Zc^*} \Big\{ \langle \zeta, X \rangle + k \int_\Omega \langle \one, \zeta(\omega) \rangle P(d\omega) - \varrho^*[\zeta] \Big\} = \sup_{\zeta \in \Zc^*} \Big\{ \langle \zeta, X \rangle - \varrho^*[\zeta] + k\varrho[\one] \Big\} \\
	        &= \varrho[X] + k\varrho[\one].
	    \end{align*}
	    Hence, the translation property is satisfied. 
	\end{proof}
	 %\Halmos\endproof

	It follows from Theorem~\ref{dual-rep} that if a risk measure $\varrho$ is lower semicontinuous and satisfies the axioms of monotonicity, convexity and translation equivariance, then representation \eqref{dual} holds with the set $\Ac$ defined as:
	\[
	\Ac = \Big\{ \zeta \in \Zc^* ~:~ \int_\Omega \zeta(\omega) dP(\omega) = \mu_\zeta, ~ \zeta \succeq 0 \Big\}. 
	\]
	\begin{corollary}
	If a risk measure $\varrho[\cdot]$ is coherent, then	
	\[
	\varrho[0] = 0 \quad\text{and}\quad \Ac= \partial \varrho[0].
	\]
	\end{corollary}
	\begin{proof}
	If $\varrho$ is also positive homogeneous, then $\varrho$ is the support function of $\Ac = \dom(\varrho^*)$. Then
	\[
	    \varrho [0] = \sup_{\zeta \in \Zc^*} \bigg\{ \langle 0,X \rangle - \varrho^*[\zeta] \bigg\} = 0.
	    \]
	   To show the form of the set $\Ac$ recall that  
	    \begin{align*} 
	    \partial \varrho[0] & = \{ \zeta \in \Zc^* : \langle \zeta, X - 0 \rangle \leq \varrho[X] - \varrho[0] \} = \{ \zeta \in \Zc^* : \langle \zeta, X \rangle \leq \varrho[X] \}. 
	    \end{align*}
	        Hence, for all $\zeta \in \Ac$, (\ref{dual}) implies that $\zeta \in \partial \varrho[0]$. On the other hand, if $\zeta \in \partial \varrho[0]$, then $\zeta \in \Ac$ by the definition of a support function.
	        \end{proof}
We shall consider further the following property.
	\begin{itemize}
	    \item[] \textbf{Normalization} a coherent measure of risk $\varrho: \Zc \to \Rb\cup\{+\infty\}$ is normalized iff $\varrho[\one] = 1$.
	\end{itemize}

	\begin{corollary}
	    For a normalized coherent measure of risk $\varrho[\cdot]$, we have $\int_\Omega \langle \onen, \zeta\rangle P(d\omega) =1$.
	\end{corollary}
	\begin{proof} 
	    It follows from equation \eqref{rhoofone} that $\langle \one, \mu_\zeta \rangle = 1$, as stated. This entails that for all $\zeta \in \Ac$, $\zeta P$ can be interpreted as a probability measure on the space $\Omega \times \{ 1, 2 ,\dots, m \}$.
	\end{proof}
	% \Halmos\endproof

	\medskip

	In the paper \cite{Burgert}, the authors have adopted the following translation axiom:

	\textbf{T.} For any constant $\alpha\in\Rb$ and any vector $e^i$ whose $i$-th component is 1 ($i=1,\dots, m$) and all other components are zero, we have $\varrho[X+\alpha e^i] = \varrho[X] +\alpha.$
	\begin{theorem}
	\label{t:comparison_Ruesch}
	    Assume that $\varrho$ is a proper lower-semicontinuous convex risk functional. Property {\rm T} holds if and only if $\int_\Omega \zeta_i(\omega) dP(\omega) =1$ for all $i=1,\dots, m$ for all $\zeta\in\dom \varrho^*$.   
	    Furthermore, if property {\rm T} holds, then 
	    \begin{gather}
	        \int_\Omega \langle \one, \zeta(\omega) \rangle dP(\omega) =\varrho[\one]=m ,\quad \int_\Omega \zeta(\omega) dP(\omega) =\onen,\quad\text{for all } \zeta\in\dom \varrho^*; \\\label{e:T-int} \\
	        \varrho[X+ a] = \varrho[X] + \varrho[a] \quad \text{for all } X\in\Zc,\; a\in\Rb^m. \\\label{e:T-randomsum}
	    \end{gather}
	\end{theorem}
	\begin{proof}
	    Suppose T holds. Then for a random vector $X$ in the domain of $\varrho$ and every $\zeta\in\Zc^*$, we have
	    \begin{align*}
	        \varrho^*[\zeta]\geq & \sup_{\alpha\in \Rb}\big\{\langle \zeta,X+\alpha e^i\rangle-\varrho[X+ \alpha e^i]\big\}
	        =\sup_{\alpha\in \Rb} \Big\{\int_\Omega \alpha\zeta_i(\omega)  P(d\omega)+\langle \zeta,X\rangle-\varrho[X] - \alpha\Big\}\\
	        & = \sup_{\alpha\in \Rb} \alpha\Big\{\int_\Omega \zeta_i(\omega)  P(d\omega)- 1\Big\}+\langle \zeta,X\rangle-\varrho[X]. 
	    \end{align*}
	    It follows that $\varrho^*[\zeta]=+\infty$ for any $\zeta\in \Zc^*$ such that
	    $\int_\Omega \zeta_i(\omega) P(d\omega) \not= 1$. 
	    This entails that for every constant vector $a\in\Rb^m$, the risk value is
	    \begin{equation}
	        \label{T2-const}
	        \varrho[a]= \varrho\Big[\sum_{i=1}^m a_ie^i\Big] = \sum_{i=1}^m a_i.
	    \end{equation}
	    The other direction is straightforward. Indeed, 
	    \begin{align*}
	        \varrho[X+\alpha e^i] &= \sup_{\zeta\in \Zc^*}\big\{\langle \zeta,X+\alpha e^i\rangle-\varrho^*[\zeta]\big\} 
	        = \sup_{\zeta\in \Zc^*} \Big\{\int_\Omega \alpha\zeta_i(\omega) P(d\omega)+\langle \zeta,X\rangle-\varrho^*[\zeta]\Big\} \\
	        & = \sup_{\zeta\in \Zc^*} \Big\{\alpha +\langle \zeta,X\rangle-\varrho^*[\zeta]\Big\}  = \varrho[X]+ \alpha.
	    \end{align*}
	    Additionally, property T also implies that 
	    \begin{equation*}
	        \int_\Omega \zeta(\omega) dP(\omega) = \onen,\quad 
	        \int_\Omega \langle \onen, \zeta(\omega) \rangle dP(\omega) = m = \varrho[\one].
	    \end{equation*}
	    Due to equation \eqref{T2-const}, for all $X\in\Zc$ and $a\in\Rb$, we obtain
	    \begin{align*}
	        \varrho[X+a] = \varrho[X+\sum_{i=1}^m a_ie^i] = \varrho[X] + \sum_{i=1}^m a_i = \varrho[X] + \varrho[a],    
	    \end{align*}
	    which completes the proof. 
	    \end{proof}
	%\Halmos\endproof
We also observe that a particular implication of Theorem~\ref{t:comparison_Ruesch} is that risk measures are linear on constant vectors.

\begin{corollary}
 	If a coherent measure of risk $\varrho[\cdot]$ satisfies property {\rm T}, then it is linear on constant vectors. 
 \end{corollary} 
\begin{proof}
	Indeed, a special case of \eqref{e:T-randomsum} shows that
	\begin{equation*}
	        \varrho[a + b] = \sum_{i=1}^m (a_i+b_i) = \varrho[a]+ \varrho[b] \quad \text{for all } a,b\in\Rb^m
	    \end{equation*}
	    This combined with the fact that $\varrho[0]=0$ and the positive homogeneity of the risk measure proves the statement.
\end{proof}

In \cite{ekeland2011law}, the authors have analyzed law-invariant risk measures for bounded random vectors. They have introduced a set of axioms that are closest to ours: their axioms include our axioms together with the two normalization properties $\varrho[\one]=1$ and $\varrho[0]=0$. We do not need these normalization properties to establish the dual representation for general random vectors with finite $p$-moments, $p\geq 1$; we derive that the risk of the deterministic zero vector is zero from the dual representation. The property of strong coherence of risk measures, introduced in that paper implies in particular that $\varrho[a + b] =  \varrho[a]+ \varrho[b],$  which appears to be a strong assumption.

% ------------------------------- LINEAR SCALARIZATION ---------------------------- %

\subsection{Risk measures obtained via sets of linear scalarizations}

	Suppose we have a random vector $X \in \Zc = \Lpm$ with a right-continuous distribution function $F(X;\cdot)$ and marginal distribution function $F_i(X_i; \cdot)$ of each component $i = 1,\dots, m$. We consider linear scalarization using vectors taken from the simplex 
	\[
	S_m^+ = \{ c \in \Rb^m ~|~ \sum_{i=1}^m c_i = 1, ~ c_i \geq 0 ~ \forall i = 1, \dots, m \}. 
	\]
	Let $\varrho:\Lc_p(\Omega,\Fc, P)\to\Rb\cup\{ +\infty\}$ be a lower semi-continuous risk measure. For any fixed set $S\subset S_m^+$, we define the risk measure 
	\begin{equation}
	    \label{rho_max_S}
	    \varrho_S[X]=\varrho[X_S],\quad\text{where } X_S (\omega)= \max_{c\in S} c^\top X(\omega),\;\; \omega\in\Omega.
	\end{equation}
	It is straightforward to see that $X_S\in\Lc_p(\Omega,\Fc, P)$ and hence, the risk measure $\varrho_S[\cdot]$ is well-defined on $\Lpm.$
	\begin{theorem}
	    If $\varrho:\Lc_p(\Omega,\Fc, P)\to\Rb\cup\{ +\infty\}$ is a coherent (convex) risk measure, then for any set $S\subset S_m^+$, the risk measure $\varrho_S[X]=\varrho[X_S]$ is coherent (convex) according to Definition~\ref{d:riskonvectors}.
	\end{theorem}
	\begin{proof}
	    For two random vectors $X,Y\in\Lpm$ with $X\leq Y$ component-wise a.s., we have $c^\top X\leq c^\top Y$ a.s. for all $c\in S_m^+$. This implies that $\max_{c\in S} c^\top X\leq \max_{c\in S} c^\top Y$ a.s. and, hence, $\varrho[X_S]\leq \varrho[Y_S]$. Thus, the monotonicity axiom is satisfied. 
	    Given two random vectors $X, Y \in \Zc $ and $\alpha \in [0,1]$, consider their convex combination $\alpha X + (1-\alpha)Y$. 
	    Due to the convexity and monotonicity of $\varrho[\cdot]$, we have
	    \begin{align*}
	        \varrho_S[\alpha X + (1-\alpha)Y] & = \varrho[\max_{c\in S} c^\top (\alpha X + (1-\alpha)Y)] \leq \varrho[\alpha \max_{c\in S} c^\top X + (1-\alpha)\max_{c\in S} c^\top Y] \\
	        & \leq \alpha \varrho[\max_{c\in S} c^\top X] + (1-\alpha)\varrho[\max_{c\in S} c^\top Y] = \alpha \varrho_S[X] + (1-\alpha)\varrho_S[Y].
	    \end{align*}
	    Thus, the convexity axiom is satisfied. Given a random vector $X \in \Zc$ and a constant $t \in \Rb$, it follows:
	    \begin{align*}
	        \varrho_S[X] & = \varrho[\max_{c\in S} c^\top (X+t\one)]  = \varrho[\max_{c\in S} (c^\top X + t c^\top \one)] = \varrho[\max_{c\in S} c^\top  X + t] = \varrho[X_S] + t
	    \end{align*}
	    Positive homogeneity follows in a straightforward manner. 
	\end{proof}
	%\Halmos\endproof

	If the set $S$ is a singleton, we obtain the following.
	\begin{corollary}
	    Let $\varrho:\Lc_p(\Omega,\Fc, P)\to\Rb\cup\{ +\infty\}$ be a coherent (convex) risk measure. For any vector $c\in\Sc_m^+$, the risk measure $\varrho_c[X]=\varrho[c^\top X]$ is coherent (convex) according to Definition~\ref{d:riskonvectors}.
	\end{corollary}

	Using the dual representation of the coherent risk measure $\varrho$ for scalar-valued random variables, we obtain the following:
	\begin{equation}
	    \begin{aligned}
	        \varrho[c^\top X] & = \sup_{\xi \in \Ac_\varrho} \int_\Omega \xi(\omega) c^\top X(\omega) P(d\omega) = \sup_{\zeta \in \At} \int_\Omega \zeta(\omega) X(\omega) P(d\omega)\quad\text{with } \At = \{ \xi c: \xi \in\Ac_\varrho\}
	    \end{aligned}
	    \label{ScalarDual}
	\end{equation}

	Additionally, a measurable selection $\nu_X(\omega) \in \argmax_{c \in S} c^\top X(\omega)$ exists by the Kuratowski-Ryll-Nadjevski theorem; we shall use the notation $\nu_X \in S$ for any such selection. 
	\begin{equation}
	    \begin{aligned}
	        \varrho[\max_{c\in S} c^\top X] &= \sup_{\xi \in \Ac_\varrho} \int_\Omega \xi(\omega) \max_{c\in S} c^\top X(\omega) P(d\omega)\notag = \sup_{\xi \in \Ac_\varrho} \int_\Omega \xi(\omega)\nu_X^\top X(\omega) P(d\omega)\notag\\
	        & = \sup_{\zeta \in \At^\prime} \int_\Omega \zeta(\omega) X(\omega) P(d\omega) \quad\text{with }\notag \At^\prime = \{ \xi \nu_X: \xi \in\Ac_\varrho\}
	    \end{aligned}
	    \label{maxScalarDual}
	\end{equation}

	Notice that the representations just derived have the form of the dual representation in \eqref{dual}, however we have not established that $\At$ coincides with the domain of its conjugate function. 

	We observe the following properties of the aggregation by a single linear scalarization.  

	\begin{proposition}
	    Given a coherent risk measure $\varrho: \Zc \to \bR$ and a scalarization vector $c\in \Sc^+_m$, for any random vector $X \in \Lpm$ risk of the vector measured by $\varrho[c^\top X]$ does not exceed the maximal risk of its components measured by $\varrho[\cdot].$ Furthermore, the following relation between aggregation methods holds
	    $ \varrho[c^\top X] \leq c^\top \varrho[X].$
	\end{proposition}
	\begin{proof}
	    The dual representation implies the following: 
	    \begin{align*} 
	        \varrho[c^\top X] &=  \sup_{\xi \in \Ac} \int_\Omega \sum_{i=1}^m c_i\xi(\omega) X_i(\omega) P(d\omega) = \sup_{\xi \in \Ac} \sum_{i=1}^m c_i \int_\Omega \xi(\omega) X_i(\omega) P(d\omega) \\
	        &  \leq  \sum_{i=1}^m \sup_{\xi \in \Ac} c_i \int_\Omega \xi(\omega) X_i(\omega) P(d\omega) = \sum_{i=1}^m c_i \varrho[X_i] \leq   \max_{1\leq i\leq m }\varrho[X_i].  
	    \end{align*}
	\end{proof}

	The penultimate relation implies the second claim of the theorem.
	% \Halmos\endproof

	We also show the following useful result, which implies that we can use statistical methods to estimate the systemic risk measure $\varrho_S[X].$ 
	\begin{proposition} 
	    \label{law-inv-max_S}
	    If $\varrho:\Lc_p(\Omega,\Fc, P)\to\Rb\cup\{ +\infty\}$ is a law-invariant risk measure, then for any set $S\subset S_m^+$, the systemic risk measure $\varrho_S[X]=\varrho[X_S]$ is law-invariant.   
	\end{proposition}
	\begin{proof} 
	    It is sufficient to show that for two random vectors $X$ and $Y$, which have the same distribution, the respective random variables $X_S$ and $Y_S$ have the same distribution. 
	\end{proof}

	We observe that $c^\top X$ and $c^\top Y$  have the same distribution for any vector $c\in\Rb^m$. Hence, for any real number $r$, the following relations hold: 
	\begin{align*}
	    P(X_S\leq r) & = P(c^\top X\leq r,\; \forall c\in S) = P(c^\top Y\leq r,\; \forall c\in S) = P(Y_S\leq r),
	\end{align*}
	which shows the equality of the distribution functions.
	% \Halmos\endproof

% ---------------------------------- NONLINEAR SCALARIZATION ----------------------------------- %

\subsection{Systemic Risk Measures Obtained via Nonlinear Scalarization}

	The second aggregation method that falls within the scope of our axiomatic framework is that of nonlinear scalarization. This class of risk measures cannot be obtained within the framework of aggregations by non-linear functions, and does not fit the axiomatic approaches in \cite{Chen} or in \cite{Burgert}. Furthermore, we shall see that this method of evaluating systemic risk allows to maintain \emph{fairness} between the system's participants.

	We define $\Omega_m = \{1,\dots,m\}$ and consider a probability space $(\Omega_m,\Fc_c, c)$, where $c\in S^+_m$ and $\Fc_c$ contains all subsets of $\Omega_m$. We view $c$ as a probability mass function of the space $\Omega_m$. Given an $m$-dimensional random vector $X \in \Zc = \Lpm$ and a collection of $m$ univariate measures of risk $\varrho_i:(\Omega,\Fc, P)\to\Rb$, $i=1,\dots, m,$  we define the random variable $X_R$ on the space $\Omega_m$ as follows:
	\begin{equation}
	    \label{XR}
	    X_R(i) = \varrho_i[X_i],\quad i=1,\dots m. 
	\end{equation}
	Choosing a scalar measure of risk $\varrho_0:(\Omega_m,\Fc_c, c)\to\Rb$, the measure of systemic risk 
	 $\varrho_s:\Lc_p(\Omega_m,\Fc_c, c)\to\Rb$ is defined as follows:
	\begin{equation}
	     \varrho_s[X] = \varrho_0[X_R]. 
	     \label{syst_risk2} 
	\end{equation} 
	This is a nonlinear aggregation of the individual risks $\varrho[X_i]$, hence this approach falls within the category of methods that evaluate the risk of each component first and then aggregate their values. The measure $\varrho_s[X]$ satisfies the axioms postulated for systemic risk measures.

	\begin{theorem}
	Suppose the univariate measures of risk $\varrho_i:(\Omega,\Fc, P)\to\Rb$, $i=1,\dots, m$ are coherent and let $\varrho_s[\cdot]$ be defined as in \eqref{syst_risk2}. Then $\varrho_s[\cdot]$ satisfies properties (A1)--(A4).
	\end{theorem}

	\begin{proof}
	    (i) Given any $X, Y \in \Zc$ and $\alpha \in (0,1)$, we consider the random vector $Z = \alpha X + (1-\alpha)Y$.  We have $\varrho_i[Z_i] \leq \alpha\varrho_i[X_i] + (1-\alpha)\varrho_i[Y_i]$, $i = 1, \dots, m$. Defining a random variable $Z^\prime$ on $\Omega_m$ by setting $Z^\prime(i) = \alpha\varrho_i[X_i] + (1-\alpha)\varrho_i[Y_i] $, we obtain that $Z_R\leq Z^\prime$. Using the monotonicity and convexity of $\varrho_0,$ we obtain
	    \[
	    \varrho_0[Z_R]\leq \varrho_0[Z^\prime] \leq \alpha\varrho_0[X_R] + (1-\alpha)\varrho_0[Y_R].
	    \]
	    Hence $\varrho_s[\alpha X + (1-\alpha)Y] \leq \alpha \varrho_s[X] + (1-\alpha)\varrho_s[Y]$. 

	    (ii) Suppose the vectors $X, Y \in \Zc$ satisfy $X\leq Y$ a.s. This implies that $X_i \leq Y_i$ a.s. and, hence, $\varrho_i[X_i] \leq \varrho_i[Y_i]$ for all $i = 1, \dots, m$ by the monotonicity property of $\varrho_i$. This further implies that $X_R\leq Y_R$, entailing that  $\varrho_0[X_R]\leq \varrho_0[Y_R]$.  Thus (A2) is satisfied.

	    (iii) Given a random vector $X \in \Zc$, $t > 0$, we have  
	    \[
	    \varrho_s[tX] = \varrho_0[(tX)_R] = \varrho_0[t (X_R)]=t\varrho_0[X_R]
	    \] 
	    where we have used the positive homogeneity property of $\varrho_i[\cdot]$ for all $i=0,1, \dots, m$. 

	    (iv) Given a random vector $X \in \Zc$ and a constant $a$, we have 
	    $(X+a\one)_R (i) = \varrho_i[X_i+a]= \varrho_i[X_i]+a$. Hence
	    $\varrho_0[(X+a\one)_R]= \varrho_0[X_R] + a.$ This shows property (A4).
	\end{proof}
	% \Halmos\endproof

	% -------------------------------------- EXAMPLES ---------------------------------------------- %
\textbf{Examples}
	\emph{A. Systemic Mean-AVaR measure}

	Consider the case when $\varrho_0$ is a convex combination of the expected value and the Average Value-at-Risk at some level $\alpha$ and all components of $X$ are evaluated by the same measure of risk $\varrho[\cdot]$. Then for any $\kappa \in [0,1]$ and $c\in S^+_m$, we have:
	\begin{align*} 
	    \varrho_s [X] & = \varrho_0[X_R] = (1-\kappa) \Eb [X_R] + \kappa \avar_\alpha [X_R]\\
	     & = (1-\kappa) \sum_{i=1}^m c_i \varrho[X_i] + \kappa \inf_{\eta \in \Rb} \Big\{  \eta + \frac{1}{\alpha} \sum_{i=1}^m c_i(\varrho[X_i] - \eta)_+ \Big\}\\
	    %& = \sum_{i=1}^m c_i \bigg((1-\kappa)\varrho[X_i]+\kappa\big(\eta +\frac{1}{\alpha} (\varrho[X_i] - \eta)_+\big)\bigg)
	\end{align*}
	Here the infimum with respect to $\eta \in \Rb$ is taken over the individual risks of the components $\varrho[X_i]$, $i = 1, \dots, m$. Hence, this method of aggregation imposes additional penalties for the components whose risk exceeds some threshold.

	\emph{B. Systemic Mean-Semideviation measure}

	Now let $\varrho_0$ be a Mean-Upper-Semideviation risk measure of the first order and all components of $X$ are evaluated by the same measure of risk $\varrho[\cdot]$. Then the measure of systemic risk can be defined as:
	\begin{align*} 
	    \varrho_s[X] & = \varrho_0[X_R] = \sum_{i=1}^m c_i \varrho[X_i] + \kappa \sum_{i=1}^m c_i \Big(\varrho[X_i] - \sum_{j=1}^m c_j \varrho[X_j] \Big)_+ \\
	    &= \sum_{i=1}^m c_i \varrho[X_i] + \kappa \sum_{i=1}^m c_i \Big( \varrho[X_i] - \sum_{j=1}^m c_j \varrho[X_j] \Big)_+
	\end{align*}
	The last representation shows that this risk measure is an aggregation of the individual risk of the components, which compares the risk of every component with the weighted average risk of all components and penalizes the deviation of the individual risk from that average.
	% \end{itemize}

	The presented method of non-linear aggregation maintains \underline{fairness} within the system and keeps the components functioning within the same level of risk.

% ---------------------------- OTHER RISK MEASURES ---------------------------------------- 

\section{Relations to multivariate measures of systemic risk}
\label{s:comparison-theoretical}

	In this section, we compare the proposed risk measures with the multivariate notions mentioned in section~\ref{s:prelim-review}. 

	Consider first the Multivariate Value-at-Risk ($\Mvar$) is given as the set of $p$-efficient points of the respective probability distribution. The following facts are shown in \cite{Dentcheva_Lai_Rusz}. For every $p \in (0,1)$ the level set $\Zp$ of a the distribution function of a random vector $X$ is nonempty and closed. For a given scalarization vector $c \geq 0$, the $p$-efficient points can be generated by solving the following optimization problem:
	\begin{equation}
	    \begin{aligned}
	        \min \quad & c^\top z \\
	        \text{s.t.} \quad & P(X \leq z) \geq p.
	    \end{aligned}
	    \label{Mvar}
	\end{equation}
	For every $c \geq 0$ the solution set of the optimization problem (\ref{Mvar}) is nonempty and contains a $p$-efficient point. Hence, given a random vector $X \in \Zc$ and a scalarization vector $c \in S_m^+$, $\Mvar$ at level $p \in (0,1)$ can be calculated as:
	\begin{align*}
	    \Mvar_{p}(X) & = \inf \big\{ c^\top  X ~|~ P(X \leq z) \geq p \big\} 
	    = \inf \big\{ c^\top X ~:~ X \in \Zp \big\}. 
	\end{align*}
	Therefore, using linear scalarizations, one can find the $p$-efficient point corresponding to any given vector $c \in S_m^+$. 

	Consider now the Multivariate Average Value-at-Risk ($\Mavar$) defined in \eqref{Prek_Mavar}. When small outcomes are preferred then, the unfavorable set of realizations of a random vector $X$ is given by the $p$-level set of $F(X;\cdot)$. Hence $\Mavar_p(X) = \Eb[\psi(X) ~|~ X \in \Zp ]$. If $X(\omega) \in \Zp$, then there exists a $p$-efficient point $v \in \Rb^m$ such that $X(\omega) \geq v$.  If the scalarization function $\psi(X)$ is monotonically nondecreasing, then $P(\psi(X) \leq \psi(v)) \geq p$. Denote the $p$-quantile of $\psi(X)$ by $\eta_X(p)$. Then we observe that 
	$ \eta_X(p) \leq \min_v \psi(v). $
	Therefore: 
	\begin{align*}
	    \Eb [\psi(X) ~|~ X \in \Zp] &= \Eb [\psi(X) ~|~ X \geq v ]
	    \geq \Eb [\psi(X) ~|~ \psi(X) \geq \eta_X(p) ] \\
	    &= \inf_\eta \bigg\{ \eta + \frac{1}{p} \Eb [(\psi(X)-\eta)_+] \bigg\}
	    = \avar_p(\psi(X)) 
	\end{align*}
	for all $p \in (0,1)$ where the cumulative distribution function of $\psi(X)$ is continuous. It follows that the Average Value-at-Risk of scalarized $X$ by a monotonically nondecreasing function $\psi(X)$ has a smaller value than $\Mavar$ defined in (\ref{Prek_Mavar}). 
	This implies in particular that for any $S\subset S_m^+$, 
	$\Mavar_p(X) \geq \varrho_S[X].$

	% --------------------------------------------------------------------------------- %

	We not turn to the Vector-valued Multivariate Average Value-at-Risk. It is calculated as one of the Pareto-efficient optimal solution of the following optimization problem:
	\begin{equation}
	    \begin{aligned}
	        \min_{\eta \in \Rb^m} \quad & \eta + \frac{1}{p} \Eb[(X-\eta)_+] \\
	        \text{s.t.} \quad & P(X \leq \eta) \geq 1-p.
	    \end{aligned}
	    \label{vmcvar}
	\end{equation}
	It is well-known that a feasible solution of a convex multiobjective optimization problem is Pareto-efficient if and only if it is an optimal solution of the scalarized problem with an objective function which is a convex combination of the multiple objectives. Then $\Vmavar$, which is the Pareto-efficient solution of the multiobjective optimization problem \ref{vmcvar}, is also optimal for the following problem:
	\begin{equation}
	    \begin{aligned}
	        \min_{v \in \Rb^m} \quad & c^\top v + \frac{1}{p} \Eb [c^\top (X-v)_+] \\
	        \text{s.t.} \quad & P(X \leq v) \geq 1-p,
	    \end{aligned}
	\end{equation}
	where $c \in \Rb^m$ is a scalarization vector taken from the simplex $S_m^+$. 

	Now for $X\in\Lpm$, we consider:
	\begin{equation*}
	    c^\top (X-v)_+  %= \sum_{i=1}^m c_i \max \{ 0, X_i - v_i \}
	      \geq \sum_{i=1}^m \max \{ 0, c_i (X_i - v_i) \} = \max \{0, c^\top (X-v) \}
	\end{equation*}
	due to the convexity of the max function. It follows that:
	\begin{equation*}
	    \inf_{v \in \Rb^m} \Big\{ c^\top v + \frac{1}{p} \Eb [c^\top (X-v)_+] \Big\} \geq \inf_{v \in \Rb^m} \Big\{ c^\top v + \frac{1}{p} \Eb [(c^\top X-c^\top v)_+] \Big\}.
	\end{equation*}

	In the scalar-valued case ($m=1$) the minimizer of the optimization problem defining $\avar_p(Z)$ is the $\var_p(Z)$ for a random variable $Z$. In the multivariate case ($m>1$), we established that the solution of (\ref{Mvar}) is the $p$-efficient point, or $\var_p(X)$, corresponding to a given scalarization vector $c \in S_m^+$. Denoting this $p$-efficient point as $v(c)$, it follows that:
	\[
	 c^\top  v(c) + \frac{1}{p} \Eb [c^\top (X-v(c))_+] \geq c^\top  v(c) + \frac{1}{p} \Eb [(c^\top X-c^\top  v(c))_+] 
	 \]
	 
	The relation $P(X \leq v(c)) \geq 1-p$ implies that $P(c^\top  X \leq c^\top  v(c)) \geq 1-p$. Denoting the $p$-quantile of $c^\top  X$ as $\eta_X(p; c)$, it follows that: $ \eta_X(p;c) \leq c^\top  v(c), $
	i.e. $\eta_X(p;c)$ is not larger than $c^\top  v(c)$. Therefore: 
	\begin{multline*}
	    \inf_{v \in \Rb^m} \Big\{ c^\top v + \frac{1}{p} \Eb[c^\top (X-v)_+] : P(X \leq v) > p \Big\} \\ \geq 
	    \inf_{v \in \Rb^m} \Big\{ c^\top v + \frac{1}{p} \Eb[(c^\top X-c^\top v)_+] \Big\} = \avar_p(c^\top X).
	\end{multline*}
	It follows that the scalarization of $\Vmavar$ results in a smaller value of the Average Value-at-Risk of the scalarized random vector, which is one of the systemic risk measures following the constructions in section \ref{s:axiomatic}.

	We do not pursue further investigation on set-valued systemic measures of risk as their calculation is numerically very expensive.

% --------------------------- NUMERICAL EXPERIMENTS ------------------------------------ %

\section{Two-stage stochastic programming problem with systemic risk measures}
\label{s:numerical}

Our goal is to address a situation, when the agents cooperate on completing a common task and risk is associated (among other factors) with the successful completion of the task. This type of situations are typical in robotics, as well as in energy systems, where the units cover the energy demand in certain area. 

\subsection{Two-stage monotropic optimization problem with a systemic risk measure}
\label{s:numerical-method}

	In this section, we consider how the proposed approaches to evaluate systemic risk can be applied to a two-stage stochastic optimization problem with a monotropic structure. Specifically, we focus on a  problem formulated as follows:
	\begin{align}
		\min_{x \in \X} ~&~ f(x) + \varrho [Q(x; \xi)] 
	\end{align}
	where $Q(x; \xi)$ has realizations $Q^s(x; \xi^s)$ defined as the optimal value of the second-stage problem in scenario $s \in S$:
	\begin{align}
		Q^s(x; \xi^s) = \min_{\mathbf{y}, z} ~&~ \sum_{i=1}^m c_i g^s_i(y_i, z) \label{p:obj} \\
		\text{s.t.} ~&~ T_i^s x + W_i^s y_i = h^s_i, ~~~ i = 1,\dots,m \label{p:dynamics} \\
		~&~ \sum_{i=1}^m A^s_i y_i = b^s \label{p:couple const} \\
		~&~ y_i \in \Y^s_i ~~~ i = 1, \dots, m \label{p:set y} \\
		~&~ B^s z \in \D^s \label{p:set z}
	\end{align}

	Here $f: \Rb^n \to \Rb$ is a continuous function that represents the cost of the first-stage decision $x \in \Rb^n$ and $\X \subset \Rb^n$ is a closed convex set. The random vector $\xi$ comprises the random data of the second-stage problem. In the second-stage problem, we would like to minimize the sum of $m$ cost functions $g_i: \Rb^l \times \Rb^p \to \Rb$ for $i = 1, \dots, m$ that depend on two second-stage decision variables: local decision variables $y_i \in \Rb^l$ for $i = 1, \dots, m$ and the common decision variable $z \in \Rb^p$. The decision variables $y_i \in \Rb^l$ are local for every $i = 1, \dots, m$, and the local constraints are represented as a closed convex set $\Y^s_i \subset \Rb^l$. The decision variable $z \in \Rb^p$ is common for all $i$ and needs global information to be calculated. The matrix $B^s$ is of size $d\times p$ and the set $\D^s \subset \Rb^d$ is a closed convex set. Note that the constraints \eqref{p:dynamics} linking the first-stage decision variable $x$ and the local second-stage decision variables $y_i$ are defined for every $i$ separately, where matrices $T^s_i \in \Rb^{k \times n}$, $W^s_i \in \Rb^{k \times l}$ and $h^s_i \in \Rb^{k}$ depend on the scenario $s$. The constraint \eqref{p:couple const} is a coupling constraint that links the local decision variables $y_i$, where $A^s_i \in \Rb^{d \times l}$ and $b^s \in \Rb^d$ depend on the scenario $s \in S$.  

	We define the total cost as the aggregation of the individual cost functions $g_i$ using some scalarization vector $c \in \Rb^m_+$ such that $\sum_{i=1}^m c_i = 1$ and we would like to develop a numerical method to solve the two-stage problem in a distributed way. Specifically, we use decomposition ideas based on the risk-averse multicut method proposed in \cite{RiskMulticut} and the multi-cut methods in risk-neutral stochastic programming to solve the two-stage problem, but we also decompose the second-stage problem into $m$ subproblems that can be solved independently in order to allow for a distributed operation of $m$ units (agents). 

	First, we discuss how to apply the decomposition method to solve the two-stage problem. We use the multicut method to construct a piecewise linear approximation of the optimal value of the second-stage problem and we approximate the measure of risk by subgradient inequalities based on the dual representation of coherent risk measures  $\varrho[Q] = \sup_{\mu \in \Ac_\varrho} \langle \mu, Q \rangle$. To this end, we introduce auxiliary variable $\eta\in\Rb$, which will contain the lower approximation of the measure of risk.  Further, we designate $Q$ the random variable with realizations $q^s$ which represent the lower approximations of the function $Q^s(\cdot,\xi^s).$
	Then the master problem in our method takes on the following form:
	\begin{equation}
	\label{p:master}
	\begin{aligned}
	    \min_{x, \eta, q} ~&~ f(x) + \eta & \\
	    \text{s.t.} ~&~ \eta \geq \langle \mu^{\tau}, Q \rangle, ~~~ \tau = 1, \dots, t-1 & \\
	    ~&~ q^{s} \geq \hat{q}^{s,\tau} + \langle g^{s,\tau}, x - x^\tau \rangle, ~~~ \tau = 1, \dots, t-1, ~ s = 1, \dots, S \\
	    ~&~ x \in \X.
	\end{aligned}
	\end{equation}
	The optimal value $\hat{\eta}^t$ contains the value of the approximation of $\varrho[Q(\hat{x}^t; \xi)],$ where $\hat{x}^t$ is the solution of the master problem at iteration $t$. Notice that the approximation $\varrho^t[Q],$ of $\varrho[Q(\hat{x}^t; \xi)],$ is given by
	\[
	\hat{\eta}^t = \varrho^t[Q] =\max_{0\leq \tau\leq t-1}  \langle Q, \mu^\tau \rangle
	\]
	with $\mu^\tau$ being the probability measures from $\Ac_\varrho$ calculated as subgradients in the previous iterations. 
	We shall explain how the subgradients $\mu^\tau$ are obtained in due course.  
	The value $\hat{q}^{s,\tau}$ is the optimal value of the second-stage problem in scenario $s$ at iteration $\tau$ and $g^{s,\tau}$ is the subgradient calculated using the optimal dual variables of the constraints $\eqref{p:dynamics}$. One can solve the second-stage problem where the objective function consists of a scalarization of $m$ cost functions, but we would like to decompose the second-stage problem into $m$ subproblems $Q^s_i$ that can be solved independently in a distributed manner.

	Consider the second-stage problem $Q^s(x;\xi^s)$ for a fixed first-stage decision variable $x \in \Rb^n$. To decompose the global problem into $m$ local subproblems, we need to handle two problems: (i) distribute the common decision variable $z \in \Rb^p$ to individual subproblems $i$; (ii) decompose the coupling constraints. The common decision variable $z$ can be distributed to subproblems by creating its copy $z_i$ for every $i$, where $i = 1, \dots, m$.
	%, with local sets $\D^s_i$ of set $\D^s$. 
	Then we ensure the uniqueness of $z$ by enforcing the decision variables $z_i$ to be equal to each other.
	Then the second-stage problem $Q(x;\xi)$ can be rewritten as:
	\begin{align}
		Q^s(x; \xi^s) = \min_{\mathbf{y}, \mathbf{z}} ~&~ \sum_{i=1}^m c_i g^s_i(y_i, z_i) \label{p2:obj} \\
		\text{s.t.} ~&~ T^s_i x + W^s_i y_i = h^s_i , ~~~ i = 1,\dots,m \label{p2:dynamics} \\
		~&~ \sum_{i=1}^m A^s_i y_i = b^s \label{p2:couple const} \\
		~&~ z_i = z_j ~~~ i,j =1,\dots,m,  \label{p2:copy const} \\
		~&~ y_i \in \Y^s_i, ~ B^sz_i \in \D^s ~~~ i = 1,\dots,m \label{p2:set y}
	\end{align}

	% Note that the constraint \eqref{p2:copy const} can be rewritten as a system of linear equations $z_i - z_j = 0$ for $i, j = 1, \dots, m$. Then one can define a matrix $M^s = [M^s_1, \dots, M^s_m]$ with appropriately chosen matrices $M^s_i$ such that $M^s_1z_1 + \dots + M^s_m z_m = 0$. 
	In order to distribute the coupling constraints \eqref{p2:couple const}, \eqref{p2:copy const}, we can apply Lagrange relaxation using Lagrange multipliers $\lambda^s \in \Rb^d$ and $\mu^s \in \Rb^{m \times m}$. Then the \textit{global augmented Lagrangian problem} $\Lambda^s_{\kappa_0}$ associated with the second-stage problem is defined as:
	\begin{align*}
		\Lambda^s_{\kappa_0} (\mathbf{y}, \mathbf{z}) &= \sum_{i=1}^m c_i g^s_i(y_i, z_i) + \langle \lambda^s, \sum_{i=1}^m A^s_i y_i - b^s \rangle + \sum_{i=1}^m \sum_{j=1}^m \mu^s_{ij}(z_i - z_j) + \frac{\kappa_0}{2} \bigg\| \sum_{i=1}^m A^s_i y_i - b^s \bigg\|^2 \\
		&+ \frac{\kappa_0}{2} \sum_{i=1}^m \sum_{\substack{j=1\\j\neq i}}^m (z_i-z_j)^2  \\
		&= \sum_{i=1}^m c_i g^s_i(y_i, z_i) + \sum_{i=1}^m \langle \lambda^s, A^s_i y_i \rangle + \sum_{i=1}^m \sum_{j=1}^m (\mu^s_{ij} - \mu^s_{ji}) z_i - \langle \lambda^s, b^s \rangle + \frac{\kappa_0}{2} \bigg\| \sum_{i=1}^m A^s_i y_i - b^s \bigg\|^2 \\
		&+ \frac{\kappa_0}{2} \sum_{i=1}^m \sum_{\substack{j=1\\j\neq i}}^m (z_i-z_j)^2
	\end{align*}
	where $\kappa_0 > 0$ is a penalty coefficient. This problem can be decomposed into subproblems with a \textit{local augmented Lagrangian} $\Lambda_{\kappa_0}^{s,i}$ defined as:
	\begin{align*}
		\Lambda^{s,i}_{\kappa_0}(y_i, \tilde{y}, z_i, \tilde{z}, \lambda, \mu) &= c_i g^s_i(y_i, z_i) + \langle \lambda^s, A^s_i y_i \rangle + \sum_{j=1}^m (\mu^s_{ij} - \mu^s_{ji}) z_i \\
		& + \frac{\kappa_0}{2} \bigg\| A^s_i y_i + \sum_{\substack{j=1\\j\neq i}}^m A^s_j \tilde{y}_j - b^s \bigg\|^2 + \frac{\kappa_0}{2} \bigg( 2 \sum_{\substack{j=1\\j\neq i}}^m (z_i - \tilde{z}_j)^2 + \sum_{k \neq i} \sum_{\substack{j=1\\j\neq k, i}}^m (\tilde{z}_k - \tilde{z}_j)^2 \bigg)
	\end{align*}
	where $(y_i, z_i)$ are decision variables of subproblem $i$, $\tilde{y}$ and $\tilde{z}$ contain given optimal decisions of other subproblems $j=1,\dots,m, ~ j \neq i$. Note that the first penalty term can be expanded as:
	\begin{align*}
		\| A^s_iy_i + \sum_{\substack{j=1\\j \neq i}}^m A_j^s\tilde{y}_j - b^s \|^2 = \sum_{k=1}^d \bigg( [A^s_iy_i]_k + \sum_{\substack{j=1\\j\neq i}}^m [A_j^s\tilde{y}_j]_k - b_k^s \bigg)^2
	\end{align*}
	This implies that for $k$ such that $[A^s_iy_i]_k = 0$, the remaining terms are constant with respect to $i$ and can be omitted from the optimization problem. Hence, the subproblem $i$ needs access to the decisions of $j = 1, \dots, m, ~ j \neq i$ which are coupled with $i$ in constraint $k$. Similarly, consider the second penalty term:
	\begin{align*}
		2 \sum_{\substack{j=1\\j\neq i}}^m (z_i - \tilde{z}_j)^2 + \sum_{k \neq i} \sum_{\substack{j=1\\j\neq k, i}}^m (\tilde{z}_k - \tilde{z}_j)^2
	\end{align*}
	The terms contained in the last summation that do not include $z_i$ are constants and can be excluded from the optimization problem. Hence, we can define the subproblem for every $i$ as follows:
	\begin{align}
		Q^s_i(x, \tilde{y}, \tilde{z}, \lambda, \mu, \xi^s) = \min_{y_i, z_i} ~&~ c_i g^s_i(y_i, z_i) + \langle \lambda, A^s_iy_i \rangle + \sum_{j=1}^m (\mu^s_{ij} - \mu^s_{ji}) z_i \\
		~&~ + \frac{\kappa_0}{2} \| A^s_iy_i + \sum_{\substack{j=1\\j \neq i}}^m A_j^s \tilde{y}_j - b^s\|^2 + \kappa_0 \sum_{\substack{j=1\\j \neq i}}^m (z_i - \tilde{z}_j)^2 \notag \\
		\text{s.t.} ~&~ T^s_ix + W^s_i y_i = h^s_i \label{sp:dynamics} \\
		~&~ y_i \in \Y^s_i,\; ~ B^sz_i \in \D^s_i \label{sp:local}
	\end{align}

	For a fixed $x$, in every scenario $s \in S$, we can implement the Accelerated Distributed Augmented Lagrangian (ADAL) method, we refer to \cite{ADAL} for detailed analysis of the method. The method consists of three main steps: (i) we solve every subproblem $Q^s_i$ to find the optimal primal variables $(\hat{y}_i^s, \hat{z}_i^s)$; (ii) update the primal variables and check if the coupling constraints \eqref{p2:couple const}, \eqref{p2:copy const} are satisfied; (iii) if the constraints are not satisfied, update the dual variables and go back to step (i). The ADAL method converges to the optimal solution $(\hat{y}^s, \hat{z}^s, \hat{\lambda}^s, \hat{\mu}^s)$ in a finitely many steps and we can calculate the optimal value of the objective function $\hat{Q}^s_i$ for every subproblem $i$. Then the global objective function of the second-stage problem can be calculated as $\hat{Q}^s(x; \xi^s) = \sum_{i=1}^m \hat{Q}^s_i(x; \xi^s) - \langle \hat{\lambda}^s, b^s \rangle$. 

	Once the second-stage problem is solved for every scenario $s$, we construct objective cuts for every scenario $s \in S$ defined as:
	$$ Q^s(x; \xi^s) \geq \hat{Q}^s(x^k; \xi^s) + \langle g^{s,k}, x - x^k \rangle $$
	where $g^{s,k}$ is the subgradient of $Q^s(x;\xi^s)$ at $x = x^k$ and scenario $s \in S$. Now note that 
	$$ \partial Q^s(x; \xi^s) = \partial \bigg[ \sum_{i=1}^m Q^s_i(x; \xi^s) - \langle \lambda, b^s \rangle \bigg] = \sum_{i=1}^m \partial Q^s_i(x; \xi^s) $$
	Hence, at $x = x^k$, the subgradient for scenario $s \in S$ can be calculated as $\partial Q^s(x^k; \xi^s) = \sum_{i=1}^m \partial Q^s_i(x^k; \xi^s)$. The subgradient $\partial Q_i^s(x^k; \xi^s)$ is given as $-(T^s_i)^\top \pi^s_i$, where $\pi^s_i$ is the Lagrange multiplier associated with the constraint \eqref{p2:dynamics} in subproblem $i$. 
	Then the proposed method for solving the two-stage problem is formulated as follows:
	\begin{itemize}
		\item[]\textbf{Step 0. } Set $t = 1$ and define initial $\mu^0 \in \Ac_\varrho$.
		\item[]\textbf{Step 1. } Solve the master problem \eqref{p:master}
			and denote its optimal solution as $(x^t, \eta^t, q^t)$.
		\item[]\textbf{Step 2. } For every scenario $s = 1, \dots, S$ apply the following method.
			\begin{itemize}
			\item[]\textit{(a)} Set $l = 1$ and define initial Lagrange multipliers $\lambda^{s,1}$, $\mu^{s,1}$ and initial primal variables $y^{s,1}, z^{s,1}$. 
			\item[]\textit{(b)} Given the Lagrange multipliers $\lambda^{s,1}, \mu^{s,1}$ and decision variables of the neighboring nodes $y^{s,l}, z^{s,l}$, every node $i$ calculates its optimal solution $(\hat{y}_i^{s,l}, \hat{z}_i^{s,l})$ by solving its local problem:
				\begin{equation}
					\label{p:local-i-s}
					\begin{aligned}
						\min_{y_i^s, z_i^s} ~&~ \Lambda^{s,i}_{\kappa_0} (y_i^s, y^{s,l}, z_i^s, z^{s,l}, \lambda^{s,l}, \mu^{s,l}) \\
						\text{s.t.} ~&~ y_i^s \in \Y^s_i, ~ B^sz_i^s \in \D^s_i 
					\end{aligned}
				\end{equation}
			\item[]\textit{(c)} Every node $i$ updates its primal variables:
				\[ y_i^{s, l+1} = y_i^{s,l} + \kappa_s (\hat{y}_i^{s,l} - y_i^{s,l}) \]
				\[ z_i^{s, l+1} = z_i^{s,l} + \kappa_s (\hat{z}_i^{s,l} - z_i^{s,l}) \]
			\item[]\textit{(d)} If the constraints 
				\begin{equation}
	        		\label{e:check-stop-s}
	        		\begin{gathered}
	        		\sum_{i=1}^m A_i^s y_i^s = \mathbf{b}^s,\quad  z_i^s = z_j^s, ~ i,j = 1,\dots,m
	         		\end{gathered}
	        	\end{equation}
	        	are satisfied, then calculate the following quantities and go to Step 3:
	        	\begin{gather*}
	        		g^{s,t} = \sum_{i=1}^m g_i^{s,l} = \sum_{i=1}^m (-T_i^s)^\top \pi_i^{s,l} \\
	        		\hat{L}^{s,t} = \sum_{i=1}^m \hat{\Lambda}^{s,i,l}_{\kappa_0} - \sum_{i=1}^m \lambda_i^{s,l} b_i^s
	        	\end{gather*}
	        	where $\pi_i^{s,l}$ is the optimal Lagrange multiplier associated with the constraint \eqref{sp:dynamics} in subproblem $i$ and $\hat{\Lambda}^{s,i,l}_{\kappa_0}$ is the optimal value of the objective function \eqref{p:local-i-s}. 

	        	If any of the constraints \eqref{e:check-stop-s} are not satisfied, update their Lagrange multipliers as follows:
	        	\begin{gather*}
	        		\lambda^{s,l+1} = \lambda^{s,l} + \kappa_0^s \kappa_s \bigg( \sum_{i=1}^m A_i^s y_i^{s,l} - b^s \bigg) \\
	        		\mu_{i,j}^{s,l+1} = \mu_{i,j}^{s,l} + \kappa_0^s \kappa_s (z_i^{s,l} - z_j^{s,l})
	        	\end{gather*}
	        	Increase $l$ by one and return to Step (b).
			\end{itemize}
		\item[]\textbf{Step 3. } Calculate $\varrho^t = \varrho[\hat{L}^t]$ and $\mu^t \in \partial \rho[\hat{L}^t]$. 
		\item[]\textbf{Step 4. } If $\varrho^t = \eta^t$, stop; otherwise, increase $t$ by one and go to Step 1. 
	\end{itemize}

	Note that the penalty parameter $\kappa^s_0$ can be chosen for every scenario $s \in S$ according to the structure of the problem. The ADAL method converges to the optimal solution in scenario $s$ if the penalty parameter $\kappa^s_0 \in (0, \frac{1}{q^s})$, where $q$ is the maximum number of nonzero rows in matrices $A^s_i$ for $i = 1, \dots, m$. Hence, $\kappa^s_0$ can be chosen close to $\frac{1}{q^s}.$

\subsection{Two-stage wireless information exchange problem}

	In this section, we formulate a two-stage information exchange problem and implement the proposed numerical method to solve it. Consider a problem in a wireless communication network consisting of $J$ robots.  We denote the team of all robots by $\mathcal{J}$. The robots collect information about the unknown environment and send the information to a set $\mathcal{K}= \{ 1, \dots, K_0 \}$ of active reporting points by multi-hop communication. The active reporting points can receive information from robots and store it.  The communication links between robots and reporting points are subject to the risk of information loss. Therefore, the objective is to choose the optimal set of active reporting points to minimize the risk associated with the amount of information lost and the proportion of the total information that was gathered but has not reached the reporting points. To this end, we shall formulate a two-stage stochastic programming problem. 

	The first-stage decision variables are known as the \emph{here and now} variables. In our problem, these are binary variables $z_k \in \{0, 1 \}$ for $k \in \mathcal{K}$, where $z_k = 1$ if the $k$-th location is selected as an active reporting point is active and $z_k = 0$ otherwise. We assume that at most $K$ reporting points can be active, where $1\leq K < K_0.$ 

	Once the reporting points are chosen, the spatial configuration of the robots is observed. We model the uncertainty of the spacial configuration and the amount of information to be observed by a set $\mathcal{S} =\{ 1,\dots S\}$ of scenarios.   
	The robots gather information about the environment and either deliver it to the reporting points or exchange it with their neighbors who then deliver it to the active reporting points. The following second-stage decision variables are involved in 
	the second-stage optimization problem. The variables $T_{ij}^s$ stand for the amount of information that is sent by node $i$ to node $j$ in scenario $s\in\mathcal{S}$. The amount of information observed but not sent by robot $i$ in scenarios $s$ is denoted by $y_i^s$. The proportion of information successfully delivered to the reporting points in scenario $s$ is denoted by $x^s$. Every robot $i$ generates information $r_i^s$ and can send it to its neighbors within some communication range. These communication links between the nodes depend on a function $R_{ij}^s$ that calculates what proportion of information sent by node $i \in \mathcal{J}$ is received and correctly decoded by node $ j \in \mathcal{J} \cup \mathcal{K}$. Then $R_{ij}^s T_{ij}^s$ is the amount of information received and correctly decoded by node $j \in \mathcal{J} \cup \mathcal{K}$ in scenario $s\in\mathcal{S}$. Then the set of neighbors of node $i$ in scenario $s \in \mathcal{S}$ can be defined as the set of nodes within its communication range $\mathcal{N}^s(i) = \{ j \in \mathcal{J} \cup \mathcal{K}: R_{ij}^s > 0 \}$.

	We associate a local risk with each robot about the information that is not communicated to neighbors or delivered to the reporting points because the information might be lost due to damage to the robot or other issues. For every $i \in \mathcal{J}$ we define:
	\[
	 y_i^s = r_i^s + \sum_{j \in \mathcal{J}} R_{ji}^s T_{ji}^s - \sum_{j \in \mathcal{J} \cup \mathcal{K}} T_{ij}^s 
	\]
	as the amount of information not communicated to neighbors nor to any of the reporting points by robot $i \in \mathcal{J}$. The systemic risk associated with the team of robots is represented by the total proportion of information not delivered to the reporting points; it is defined as $(1-x^s)$, where $x^s$ is calculated as follows:
	\begin{equation}
	\label{e:proportion-centralized}
	x^s \sum_{i \in \mathcal{J}} r_i^s = \sum_{i \in \mathcal{J}}\sum_{k \in \mathcal{K}} R_{ik}^s T_{ik}^s.
	\end{equation}
	To implement the distributed method for the operation of the robots, we introduce copies of the total proportion variable for each robot (denoted $x_i^s$). We then introduce additional constraints to impose equality among the auxiliary variables $x_i^s$. Constraint \eqref{e:proportion-centralized} is then replaced by the following set of constraints:
	\begin{gather*}
	\sum_{i \in \mathcal{J}} x_i^s r_i^s = \sum_{j \in \mathcal{J}} \sum_{k \in \mathcal{K}} R_{jk}^s T_{jk}^s \\
	 x_i^s = x_j^s \quad \forall i, j \in \mathcal{J}.
	 \end{gather*}
	Using these variables, we can express the loss function of every robot $i \in \mathcal{J}$ in scenario $s \in \mathcal{S}$ as follows
	\[
	q_i^s = c_1 y_i^s + c_2 (1-x_i^s),
	\]  
	where $c_1>0$ is the weight associated with the local risk, while $c_2>0$ is the systemic risk. These are modeling parameters. We have used a choice of $c_1+c_2=1$ in-line with our theoretical proposal for aggregating sources of risk.
	The first-stage optimization problem takes on the following form:
	\begin{align*}
	    \min_{z} ~&~ \varrho[Q(z; \xi)] \\
	    \text{s.t.} ~&~ \sum_{k \in \mathcal{K}} z_k \leq K \\
	    ~&~ z_k \in \{ 0, 1 \}, ~~~ k \in \mathcal{K}.
	\end{align*}
	Here $Q(z; \xi)$ is a random variable with realizations $Q^s(z; \xi^s)$ for $s = 1, \dots, S$ denoting the optimal values of the second-stage problem. The second-stage problem deals with the operation of the robots after the location of the reporting points is fixed; it is formulated as follows:
	\begin{align}
	Q^s(z; \xi^s) &~ = \notag\\
	    \min_{T^s, y^s, x^s} &~ \sum_{i \in \mathcal{J}} w_i q_i^s \quad \text{s.t.} \label{p:obj1} \\
	    ~&~ y_i^s = r_i^s + \sum_{j \in \mathcal{J}} R_{ji}^s T_{ji}^s - \sum_{j \in \mathcal{J} \cup \mathcal{K}} T_{ij}^s, ~~~ i \in \mathcal{J} \label{p:flow}\\
	    ~&~ \sum_{j \in \mathcal{J}} T_{ji}^s + \sum_{j \in \mathcal{J} \cup \mathcal{K}} T_{ij}^s \leq a, ~~~ i \in \mathcal{J}
	    \label{p:capacity} \\
	    ~&~ \sum_{i \in \mathcal{J}} x_i^s r_i^s = \sum_{j \in \mathcal{J}} \sum_{k \in \mathcal{K}} R_{jk}^s T_{jk}^s \label{p:proportion}\\
	    ~&~ x_i^s = x_j^s, ~~~ i, j \in \mathcal{J} \label{p:consistency}\\
	    ~&~ T_{ik}^s \leq M z_k, ~~~ i \in \mathcal{J}, k \in \mathcal{K} \label{p:logicalUB}\\
	    ~&~ T_{ij}^s \geq 0, ~~~ i \in \mathcal{J}, j \in \mathcal{J} \cup \mathcal{K} \\%\label{p:nonneg-T}\\
	    ~&~ y_i^s \geq 0, ~~~ i \in \mathcal{J}.%\label{p:nonneg-y}
	\end{align}
	Here $M>0$ is a large constant, which helps us provide a logical upper bound for communication only to the active reporting points in constraint \eqref{p:logicalUB}.  Additionally, $w_i>0$ are weights associated with the loss functions of individual robots. Here again, we propose $\sum_{i\in\mathcal{J}} w_i=1$ according to the proposed systemic measures of risk. 
	We notice that the second-stage problems are always feasible for any feasible first-stage decision. Hence, the two-stage problem has a relatively complete recourse. Furthermore,  
	the recourse function $Q(z; \xi) $ has finitely many realizations $Q^s(z; \xi^s),$ $s = 1, \dots, S$ for every fixed argument $z$. This implies that we can use any coherent measure of risk $\varrho[\cdot]$ for evaluating the risk of $Q(\cdot,\xi)$; the value $\varrho[Q(z; \xi)]$ is well defined and finite for all feasible first-stage decisions $z.$ 

\subsection{Numerical results}

	We solve the problem using the distributed method proposed in \ref{s:numerical-method}. In the given problem, the decision variables $T_{ij}^s, y_i^s, x_i^s$ are local to every node $i \in \mathcal{J}$ and there are four coupling constraints that need to be distributed to the nodes: 
	\begin{itemize}
	    \item[(1)] the flow conservation constraints \eqref{p:flow};
	    \item[(2)] the transmission capacity constraints \eqref{p:capacity};
	    \item[(3)] the proportion constraint \eqref{p:proportion};
	    \item[(4)] the equality constraints enforcing the uniqueness of the proportion \eqref{p:consistency}.
	\end{itemize}

	Since the ADAL algorithm operates equality constraints, we can introduce auxiliary variables $u_i^s$ and redefine \eqref{p:capacity} as $\sum_{j \in \mathcal{J}}T_{ji}^s s+ \sum_{j \in \mathcal{J} \cup \mathcal{K}} T_{ij}^s + u_i^s = a $. Then the coupling constraints can be rewritten using appropriate matrices and vectors stacking decision variables of every node $i \in \mathcal{J}$. Let $v_i^s = [y_i^s, T_{i1}^s, T_{i2}^s, T_{i3}^s, \dots, T_{iJ}^s, \dots, T_{i(J+K_0)}^s, x_i^s, u_i^s]^\top$ be a $(J+K_0+3)$-dimensional vector stacking all decision variables of the node $i \in \mathcal{J}$. Then the flow conservation constraint \eqref{p:flow} can be rewritten as $\sum_{i \in \mathcal{J}} A_i^s v_i^s = \mathbf{r}^s$, where $\mathbf{r}^s = [r_1^s, \dots, r_J^s]^\top$ is a $J$-dimensional vector and $A_i^s$ is $J \times (J + K_0 + 3)$ dimensional matrix defined as:
	$$ A_i^s = \begin{bmatrix}
		1 & 1-R_{ii}^s & 1 & 1 & \dots & 1 & \dots & 1 & 0 & 0 \\
	    0 & 0 & -R_{i2}^s & 0 & \dots & 0 & \dots & 0 & 0 & 0 \\
	    0 & 0 & 0 & -R_{i3}^s &\dots & 0 & \dots & 0 & 0 & 0 \\
	    \vdots & \vdots & \vdots &  & \vdots &  & \vdots & \vdots \\
	    0 & 0 & 0 & 0 & \dots & -R_{iJ}^s & \dots & 0 & 0 & 0\\
	\end{bmatrix} $$
	where all elements in the $i$-th row are equal to 1 except of terms $(i,i+1)$ and the last two columns which are equal to $0$. Note that $R_{ii}^s = 1$ for all $i \in \mathcal{J}$, hence the terms $1-R_{ii}^s$ are equal to $0$. Similarly, constraints \eqref{p:capacity} and \eqref{p:proportion} can be rewritten as 
	\begin{gather*} 
	\sum_{i \in \mathcal{J}} B_i^s v_i^s = a \mathbf{1},\quad
	\sum_{i \in \mathcal{J}} (C_i^s)^\top v_i^s = 0, \quad
	\end{gather*}
	using appropriate matrices $B_i^s, C_i^s$ for every node $i \in \mathcal{J}$, where $\mathbf{1}$ is a vector of all ones. Note that since nodes can share information only with the neighbors, one can enforce the equality of the proportion variable between neighboring nodes and rewrite the constraint \eqref{p:consistency} as follows:
	\begin{align}
		x_i^s = x_j^s, ~~~ i \in \mathcal{J}, ~ j \in \mathcal{N}^s(i) \label{p:new consistency}
	\end{align}
	where $\mathcal{N}^s(i)$ is the set of nodes within communication range of node $i \in \mathcal{J}$ in scenario $s \in \mathcal{S}$. If the network is connected, constraint \eqref{p:new consistency} enforces all $x_i$ to be equal to each other and ensures the consistency and uniqueness of $x_i$. 

	We assume that the team of robots works on a square map given by the points with relative coordinates $(0,0)$ and $(2,2)$. The spatial distribution of available information to be gathered follows a normal distribution with an expected value $\mathcal{C} = (0.5, 1.75)$ in the upper left corner of the map. The network consists of 50 robots and 4 potential locations of the reporting points. We generated 200 scenarios for different spatial configurations of the robots. The four potential locations for the reporting points are fixed in the positions $(0.5,0.3), (1.5,0.25), (1.75, 0.5), (1,0.2)$. The rate function $R_{ij}^s$ depends on the distance between the nodes in the network and is defined as:
	\[ R_{ij}^s =
	    \begin{cases}
	        1, & \text{if } \| d_{ij}^s \| \leq \ell, \\
	        a \| d_{ij}^s \|^3 + b\|d_{ij}^s\|^2 + c\|d_{ij}^s\| + e, & \text{if } \ell < \|d_{ij}^s\|  \leq u, \\
	        0, & \text{if } \|d_{ij}^s\| > u,
	    \end{cases}
	\]
	where $\|d_{ij}^s\|$ is the distance between the nodes in scenario $s \in \mathcal{S}$. We set $\ell = 0.3$ and $u = 0.6$, and values $a, b, c$ and $e$ are chosen so that $R_{ij}^s$ is a continuous function. This function is commonly used in literature, see e.g. \cite{netintegrity}. The information $r_i^s$ gathered by robot $i$ in scenario $s$, depends on the robot's position relative to the expected value $\mathcal{C}$ given above. In our experiments $r_i^s$ is calculated as follows:
	\[
	r_i^s = \frac{w}{2\pi | \Sigma |^\frac{1}{2}} e^{-\frac{1}{2}(d_i - \mathcal{C}) \Sigma^{-1}(d_i - \mathcal{C})}, 
	\]
	where $d_i$ is the positions of robot $i \in \mathcal{J},$ $w$ is a scaling factor, and  $\Sigma$ is a covariance matrix, which keep fixed for all experiments.

	\begin{figure*}[!htb]
	    \centering
	    \begin{subfigure}[t]{0.3\textwidth}
	        \centering
	        \captionsetup{width=1\linewidth}
	        \includegraphics[width=\textwidth]{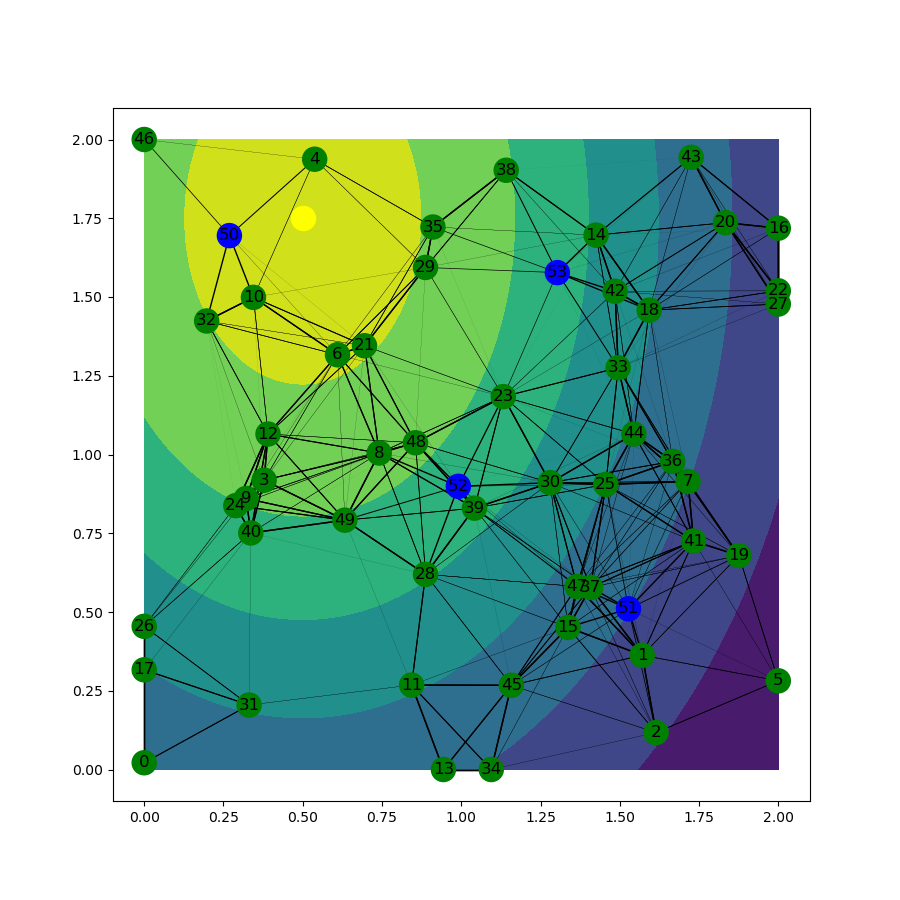}
	        \subcaption{}
	    \end{subfigure}
	    \hspace{1em}
	    \begin{subfigure}[t]{0.3\textwidth}
	        \centering
	        \captionsetup{width=1\linewidth}
	        \includegraphics[width=\textwidth]{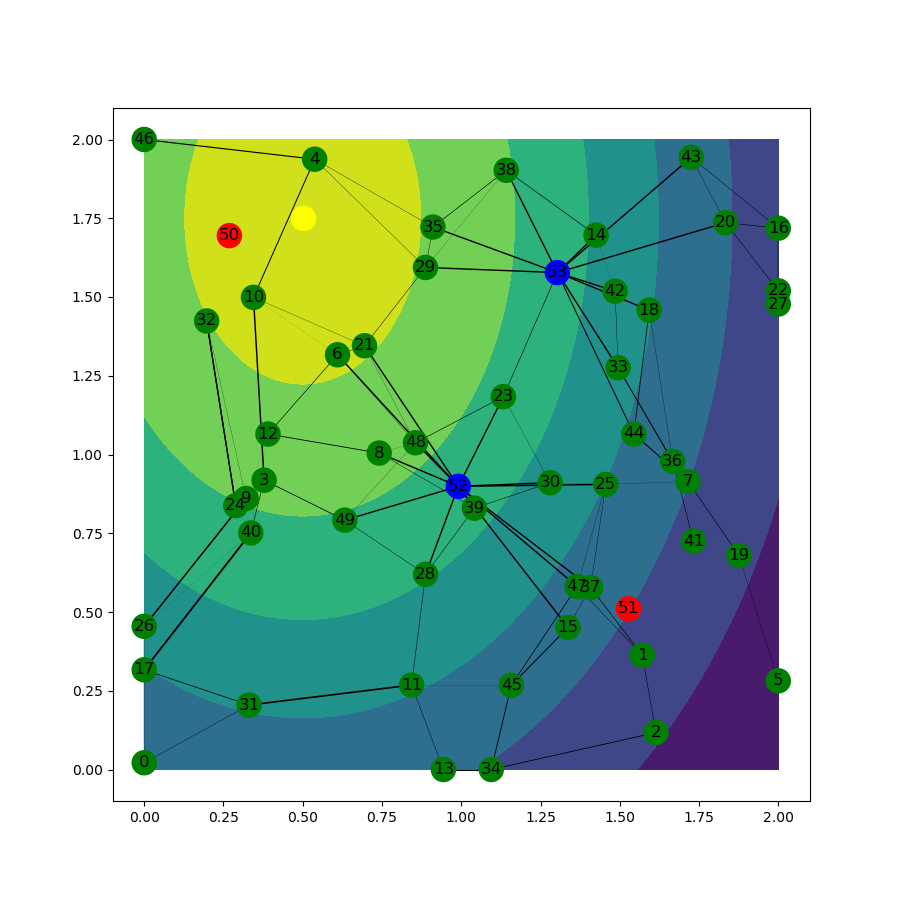}
	        \subcaption{}
	    \end{subfigure}
	    \hspace{1em}
	    \begin{subfigure}[t]{0.3\textwidth}
	        \centering
	        \captionsetup{width = 1\linewidth}
	        \includegraphics[width=\textwidth]{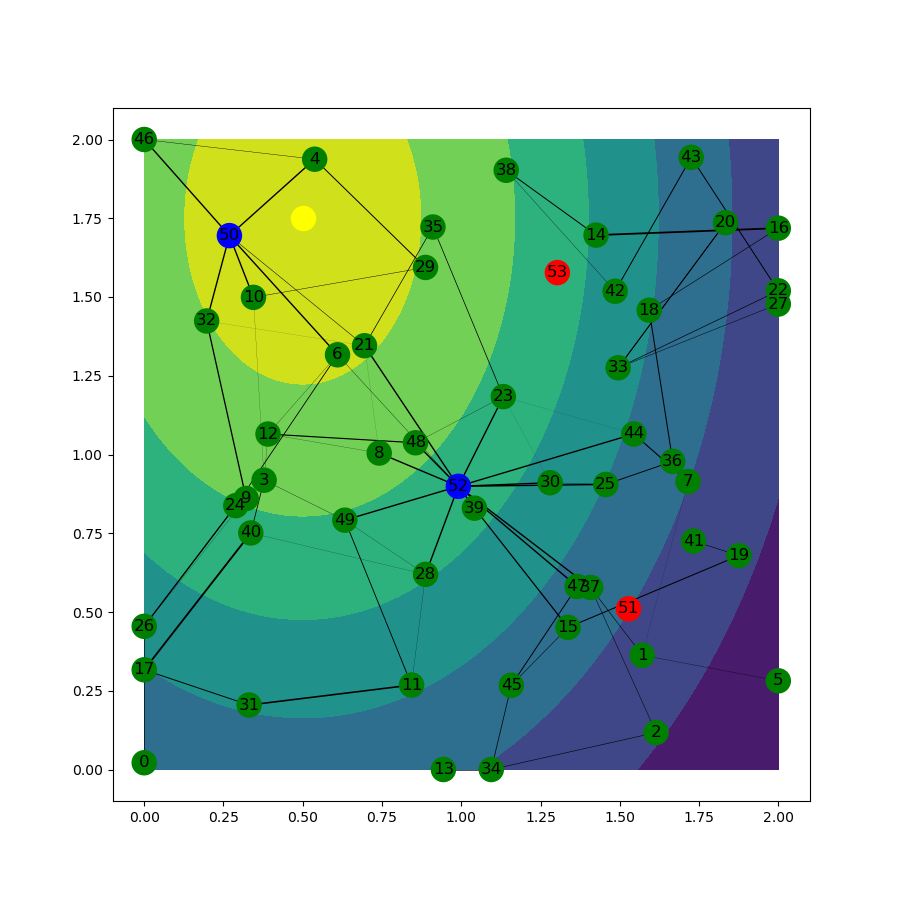}
	        \subcaption{}
	    \end{subfigure}
	    \caption{Communication network of 50 robots and 4 reporting points in one scenario. The source is located in the upper left corner. The lighter color (yellow) and darker color (purple) indicate higher and lower rates of information generation, respectively. (a) The initial spatial configuration of robots (green) and the reporting points (blue). The lines show communication links and their thickness indicates the rate of connection $R_{ij}$ between nodes $i$ and $j$. (b) The optimal routing decisions between nodes when the risk is applied to the total loss of all nodes. Blue nodes are selected and red nodes are not selected. (c) The optimal routing decisions between nodes when the individual risks of the nodes are aggregated. Blue nodes are selected and red nodes are not selected.}
	    \label{fig:info_exch}
	\end{figure*}

	\textbf{Comparison of aggregation methods.} We solved the optimization problem using two different aggregation methods:
	\begin{itemize}
	\item \emph{aggregate first}  Using the proposed multivariate measures of risk, we
	 aggregate the individual losses of the robots with a fixed scalarization $w$, we calculate for each scenario 
	 $V^s = \sum_{i \in \mathcal{J}} w_i q_i^s$  and evaluate its risk $\varrho[V]$ by several scalar-valued measures of risk; 
	 \item \emph{evaluate first} We evaluate the individual risk of every robot across all scenarios and calculate $V_i = \varrho_i[q_i]$. Then we aggregate their values $\varrho_S[V]$ using two examples of nonlinear aggregation shown in section 4.2.1.
	\end{itemize}
	  We solve the problem using a linear scalarization vector $w$ with equal weights $w_i = \frac{1}{J}$ for all $i \in \mathcal{J}$, $c = [0.8, 0.2]$ and $\avar_\alpha(\cdot)$ for three values of $\alpha = 0.1, 0.2, 0.3$. 

	\begin{figure*}[!ht]
	    \centering
	    \begin{subfigure}[t]{0.3\textwidth}
	        \centering
	        \captionsetup{width=1\linewidth}
	        \includegraphics[width=\textwidth]{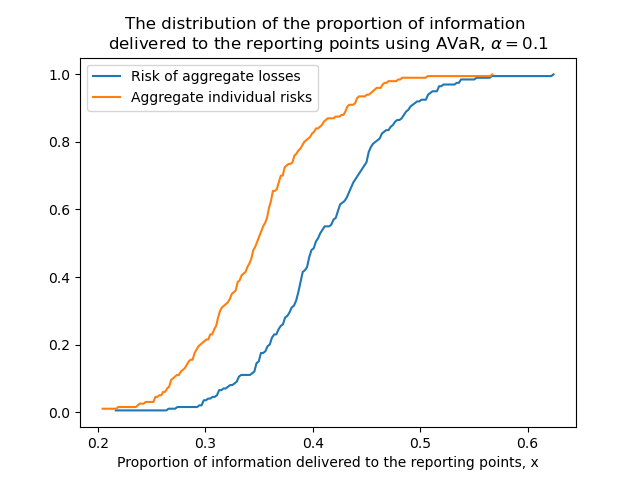}
	        \subcaption{}
	    \end{subfigure}
	    \hspace{1em}
	    \begin{subfigure}[t]{0.3\textwidth}
	        \centering
	        \captionsetup{width=1\linewidth}
	        \includegraphics[width=\textwidth]{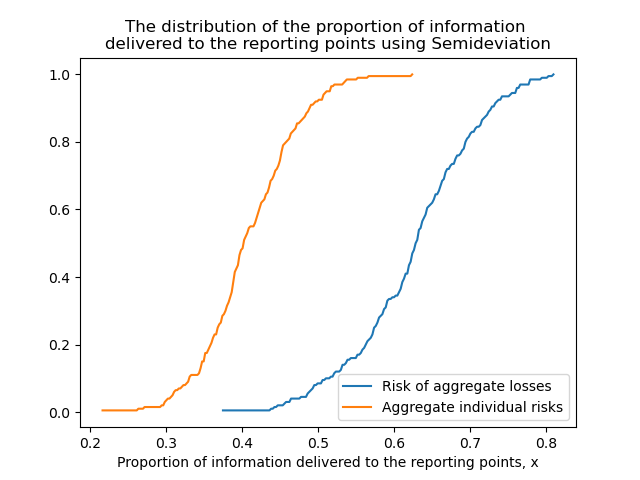}
	        \subcaption{}
	    \end{subfigure}
	    \hspace{1em}
	    \begin{subfigure}[t]{0.3\textwidth}
	        \centering
	        \captionsetup{width=1\linewidth}
	        \includegraphics[width=\textwidth]{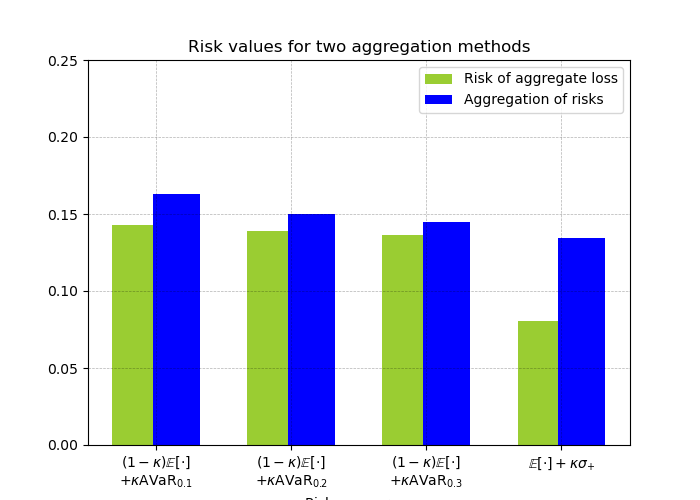}
	        \subcaption{}
	    \end{subfigure}
	    \caption{(a) Proportion of information delivered to the reporting points using Mean-AVaR at $\alpha = 0.1$. (b) Proportion of information delivered to the reporting points using Mean-Upper-Semideviation of order 1. (c) Comparison of the risk values for two aggregation methods.}
	    \label{fig:proportion}
	\end{figure*}

	The setup of the communication network problem and the optimal solutions in one of the scenarios for both methods are shown in Fig. \ref{fig:info_exch}. One can notice that depending on what kind of aggregation method is used, the set of optimal reporting points might be different. The distribution of the proportion $x$ of information delivered to the reporting points for two methods is shown in Fig. \ref{fig:proportion}. It can be seen that more information is delivered to the reporting points if we aggregate the losses of robots and evaluate the risk. This observation is also reflected in the values of the risk for both methods: imposing a risk measure on linear scalarization of the individual losses results in smaller values than aggregation of individual risks. 

	Using the optimal values of the decision variables, we can calculate $\Mavar$ and $\Vmavar$ to compare their values with the AVaR applied on linear scalarization of the random cost. The following formulas were used to calculate the values:
	\begin{gather}
	    \avar_\alpha(V) = \inf_{\eta \in \Rb} \bigg\{ \eta + \frac{1}{\alpha} \Eb \Big[ (w^\top q - \eta)_+ \Big] \bigg\} 
	    \label{num_avar}\\
	    \Mavar_\alpha(V) = \Eb \Big[ w^\top  q ~|~ q \in \mathcal{Z}_{1-\alpha} \Big]
	    \label{num_mavar}\\
	    \Vmavar_\alpha(V) = \inf_{\eta \in \Rb^2} \bigg\{ w^\top  \eta + \frac{1}{\alpha} \Eb \Big[ w^\top  (V - \eta)_+ \Big] : \Pr[V \leq \eta] > 1-\alpha \bigg\}
	    \label{num_vmcvar}
	\end{gather}
	 The values of $\avar$, $\Mavar$ and $\Vmavar$ are shown in Table \ref{table: avar_vals}. It can be seen that $\avar_\alpha(V)$ results in smaller values than $\Mavar_\alpha(V) $ and $\Vmavar_\alpha(V)$ at all confidence levels $\alpha$ as it was shown theoretically in section \ref{s:comparison-theoretical}.Those measures of risk are computationally very demanding and not amenable to the type of decision problems, we are considering. Hence, we only compare their values for the decision obtained via our proposed method. 

	\begin{table}[]
	    \centering
	    \begin{tabular}{ |c|c|c|c| } 
	    \hline
	    \diagbox{$\varrho$}{$\alpha$} & 0.1 & 0.2 & 0.3 \\ \hline
	    $\avar_\alpha$ & 0.1429 & 0.1389 & 0.1364 \\ 
	    $\Mavar_\alpha$ & 0.1992 & 0.1693 & 0.1634 \\  
	    $\Vmavar_\alpha$ & 0.174 & 0.1622 & 0.1553 \\ 
	    \hline
	    \end{tabular}
	    \caption{Comparison of $\avar$, $\Mavar$ and $\Vmavar$ values for $\alpha = 0.1, 0.2, 0.3$.}
	    \label{table: avar_vals}
	\end{table}

\begin{figure}[!ht]
	    \centering
	    \begin{subfigure}[t]{0.4\textwidth}
	        \centering
	        \captionsetup{width=1\linewidth}
	        \includegraphics[width=\textwidth]{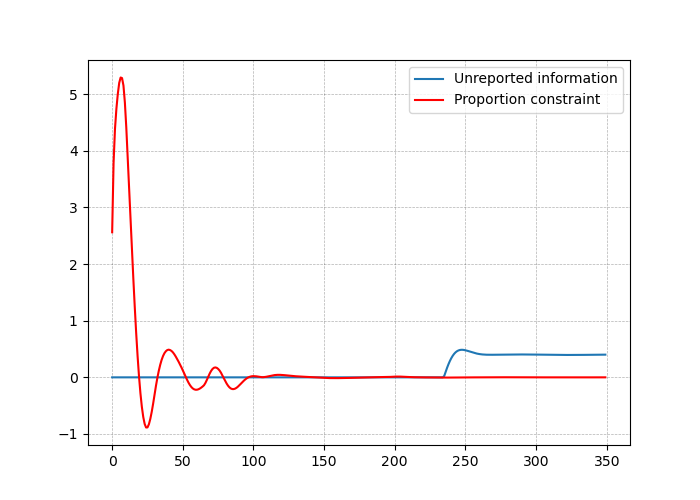}
	        \subcaption{}
	    \end{subfigure}
	    \hspace{1em}
	    \begin{subfigure}[t]{0.4\textwidth}
	        \centering
	        \captionsetup{width=1\linewidth}
	        \includegraphics[width=\textwidth]{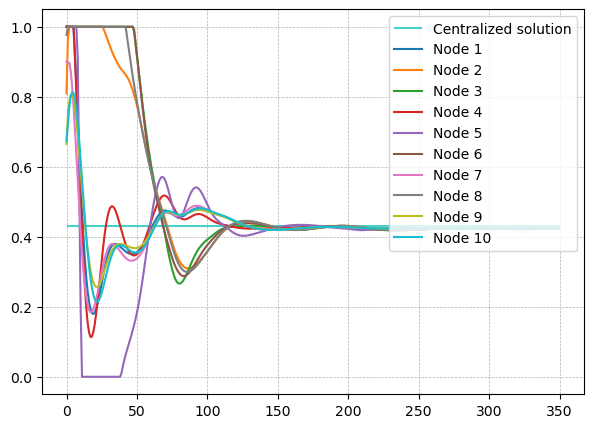}
	        \subcaption{}
	    \end{subfigure}
	    \caption{(a) Evolution of the sum of local losses $\sum_{i \in \mathcal{J}} y_i^s$ (blue) and proportion constraint $\sum_{i \in \mathcal{J}} \Big( x_i^s r_i^s - \sum_{k \in \mathcal{K}} R_{ik}^s T_{ik}^s \Big)$ (red) in scenario $s$. (b) Convergence of the robots' proportion variables $x_i^s$ to the centralized solution in scenario $s$.}
	    \label{fig:distributed}
	\end{figure}

	When we solve the problem in a distributed way, we use a smaller network consisting of 20 robots and 4 reporting points in a 1.5 by 1.5 square over 100 scenarios. It is assumed that the network is connected in all possible scenarios, that is, every node has at least one neighbor within the communication range, and all nodes are connected to the reporting points through multiple hops. This assumption is necessary for the proper calculation of the proportion of information delivered to the reporting points. If one of the nodes is isolated from the network, the rest of the group converges to a solution that does not take into account the isolated node's contribution. The problem is solved in both centralized and distributed ways, and the results for one of the scenarios are shown in Fig. \ref{fig:distributed}. As it can be seen in Fig. \ref{fig:distributed} (b), nodes converge to the centralized solution of the proportion of information delivered to the reporting points. 

\section{Conclusions}
\label{s:conclusions}

	Our contributions can be summarized as follows. 
	% \begin{itemize}
	% \item[1.] 
	We propose a sound axiomatic approach to measures of risk for distributed systems. We show that several classes of non-trivial measures that satisfy the axioms can be constructed. These measures can be calculated efficiently and are less conservative than most of the other systemic measures of risk. The class of measures proposed in section~3.3 goes beyond the popular ways to evaluate risk of agents then aggregate.

	We have devised a distributed method for solving the risk-averse two-stage problems with monotropic structure, which works for any measure of risk not only for those that are representable as expected value. The construction is quite general and could serve as a template for devising other distributed methods for problems with systemic measures of risk. 
	
	% \item[2.] 
	We demonstrate the viability of the proposed framework on a non-trivial two-stage problem involving wireless communication. The numerical experiments confirm the theoretical observations and show the advantage of the proposed approach to risk aggregation in distributed systems.

	% \item[3] 
	
	% \end{itemize}

	In conclusion, the advantage of the new approach is the good balance of robustness to the uncertainty, optimality of the loss functions involved, and the efficiency of the numerical operation.


\begin{thebibliography}{99}
\bibitem{Ararat}
Çağın Ararat, and Birgit Rudloff.
\newblock Dual representations for systemic risk measures.
\newblock {\em Mathematics and Financial Economics}, 14:139-174, 2020.

\bibitem{Artzner}
Philippe Artzner, Freddy Delbaen, Jean-Marc Eber, and David Heath.
\newblock Coherent measures of risk.
\newblock {\em Mathematical Finance}, 9(3):203--228, July 1999.

\bibitem{Biagini}
Francesca Biagini, Jean-Pierre Fouque, Marco Frittelli, and Thilo
  Meyer-Brandis.
\newblock \href{https://arxiv.org/abs/1503.06354}{A Unified Approach to
  Systemic Risk Measures via Acceptance Sets}.
\newblock 2015.

\bibitem{Brunnermeier}
Markus Brunnermeier and Patrick Cheridito.
\newblock \href{https://www.mdpi.com/2227-9091/7/2/46}{Measuring and Allocating
  Systemic Risk}.
\newblock {\em Risks}, 7(2), 2019.

\bibitem{Burgert}
Christian Burgert and Ludger Rüschendorf.
\newblock Consistent risk measures for portfolio vectors.
\newblock {\em Insurance: Mathematics and Economics}, 38(2):289--297, 2006.

\bibitem{ADAL}
Nikolaos Chatzipanagiotis, Darinka Dentcheva, and Michael Zavlanos.
\newblock An augmented lagrangian method for distributed optimization.
\newblock {\em Mathematical Programming}, 152, 01 2014.

\bibitem{Chen}
Chen Chen, Garud Iyengar, and C.~Ciamac Moallemi.
\newblock
  \href{https://pubsonline.informs.org/doi/epdf/10.1287/mnsc.1120.1631}{An
  Axiomatic Approach to Systemic Risk}.
\newblock {\em Management Science}, 59(6):1373--1388, 2013.

\bibitem{Delbaen}
Freddy Delbaen.
\newblock Coherent risk measures on general probability spaces.
\newblock {\em Advances in Finance and Stochastics}, pages 1--37, March 2000.

\bibitem{Dentcheva_Lai_Rusz}
Darinka Dentcheva, Bogumila Lai, and Andrzej Ruszczy\'nski.
\newblock Dual methods for probabilistic optimization.
\newblock {\em Mathematical Methods of Operations Research}, 60:331--346, 2004.

\bibitem{ekeland2011law}
Ivar Ekeland  and Walter Schachermayer, Law invariant risk measures on $\Lc^\infty (\Rb^d)$
\emph{Statistics \& Risk Modeling},
{28}({3}): {195--225}, {2011}.

\bibitem{Kromer}
Kromer Eduard, Overbeck Ludger, and Zilch Katrin.
\newblock \href{https://dx.doi.org/10.2139/ssrn.2268105}{Systemic Risk Measures
  on General Measurable Spaces}.
\newblock 05 2016.

\bibitem{Feinstein}
Zachary Feinstein, Birgit Rudloff, and Stefan Weber.
\newblock \href{https://arxiv.org/abs/1502.07961}{Measures of Systemic Risk}.
\newblock {\em {SIAM} Journal on Financial Mathematics}, 8(1):672--708, Jan
  2017.

\bibitem{FollmerSchied}
Hans F\"{o}llmer and Alexander Schied.
\newblock Convex measures of risk and trading constraints.
\newblock {\em Finance and Stochastics}, 6:429--447, 2002.

\bibitem{Follmer}
Hans F\"{o}llmer and Alexander Schied.
\newblock {\em Stochastic Finance: An Introduction in Discrete Time, 3rd
  Edition}.
\newblock Walter De Gruyter, 2011.

\bibitem{RiskMulticut}
Sıtkı G\"ulten and Andrzej Ruszczy\'{n}ski.
\newblock \href{https://doi.org/10.1007/s10479-014-1768-2}{Two-stage portfolio
  optimization with higher-order conditional measures of risk}.
\newblock {\em Annals of Operations Research}, 229:409--427, 06 2015.

\bibitem{hamel2010duality}
Andreas Hamel, and Frank Heyde, {Duality for set-valued measures of risk},
\emph{SIAM Journal on Financial Mathematics},
{1}({1}):{66--95},{2010}.

\bibitem{jouini2004vector}
Elyes Jouini, Moncef Meddeb, and  Nizar Touzi, 
{Vector-valued coherent risk measures},
\emph{Finance and stochastics},
 {8}(4): {531--552}, {2004}.

\bibitem{KijOhn:1993}
M.~Kijima and M.~Ohnishi.
\newblock Mean-risk analysis of risk aversion and wealth effects on optimal
  portfolios with multiple investment opportunities.
\newblock {\em Ann. Oper. Res.}, 45:147--163, 1993.

\bibitem{LeePrekopa}
Jinwook Lee and Andr{\'a}s Pr{\'e}kopa.
\newblock
  \href{https://link.springer.com/article/10.1007/s10479-013-1482-5}{Properties
  and calculation of multivariate risk measures: MVaR and MCVaR}.
\newblock {\em Annals of Operations Research}, 211:225--254, 2013.

\bibitem{Leitner:2005}
J.~Leitner.
\newblock A short note on second-order stochastic dominance preserving coherent
  risk measures.
\newblock {\em Mathematical Finance}, 15:649--651, 2005.

\bibitem{Ma-me-Zavlanos}
Ma, Wann-Jiun and Oh, Chanwook and Liu, Yang and Dentcheva, Darinka and Zavlanos, Michael M.,
\newblock Risk-Averse Access Point Selection in Wireless Communication Networks, 
\newblock \emph{IEEE Transactions on Control of Network Systems}, 6(1): 24--36, 2019.

\bibitem{Merakli}
Merve Merakl{\i} and Simge Kü{\c{c}}ükyavuz.
\newblock \href{https://arxiv.org/abs/1708.01324}{Vector-Valued Multivariate
  Conditional Value-at-Risk}.
\newblock {\em Operations Research Letters}, 46(3):300--305, 2018.

\bibitem{Noyan2013optimization}
Nilay Noyan and G{\'a}bor Rudolf.
\newblock Optimization with multivariate conditional value-at-risk constraints.
\newblock {\em Operations research}, 61(4):990--1013, 2013.

\bibitem{pflug2018systemic}
Georg Pflug and Alois Pichler, {Systemic risk and copula models},
\emph{Central European Journal of Operations Research},
{26}: {465--483} 2018.

\bibitem{PflRom:07}
G.Ch. Pflug and W.~R\"omisch.
\newblock {\em Modeling, Measuring and Managing Risk}.
\newblock World Scientific, Singapore, 2007.

\bibitem{Prekopa}
Andr{\'a}s Pr{\'e}kopa.
\newblock
  \href{https://link.springer.com/article/10.1007/s10479-010-0790-2#citeas}{Multivariate
  Value at Risk and Related Topics}.
\newblock {\em Annals of Operations Research}, 193:49--69, 2012.

\bibitem{rockafellar2013fundamental}
R~Tyrrell Rockafellar and Stan Uryasev.
\newblock The fundamental risk quadrangle in risk management, optimization and
  statistical estimation.
\newblock {\em Surveys in Operations Research and Management Science},
  18(1-2):33--53, 2013.

  \bibitem{RuschendorfBook}
Ludger R\"{u}schendorf.
\newblock {\em Mathematical Risk Analysis.}
\newblock Springer, 2015. 

\bibitem{RS:2005}
A.~Ruszczy\'{n}ski and A.~Shapiro.
\newblock Optimization of risk measures.
\newblock In G.~Calafiore and F.~Dabbene, editors, {\em Probabilistic and
  Randomized Methods for Design under Uncertainty}, pp. 117--158,
  Springer-Verlag, London, 2005.

\bibitem{RuSh:2006a}
A.~Ruszczy\'{n}ski and A.~Shapiro.
\newblock Optimization of convex risk functions.
\newblock {\em Mathematics of Operations Research}, 31:433--452, 2006.

\bibitem{textbook}
Alexander Shapiro, Darinka Dentcheva, and Andrzej Ruszczy\'{n}ski.
\newblock {\em
  \href{https://epubs.siam.org/doi/book/10.1137/1.9780898718751}{Lectures on
  Stochastic Programming: Modeling and Theory}}.
\newblock Society for Industrial {\&} Applied Mathematics ({SIAM}), 2009.

\bibitem{netintegrity}
Michael~M. Zavlanos, Alejandro Ribeiro, and George~J. Pappas.
\newblock Network integrity in mobile robotic networks.
\newblock {\em IEEE Transactions on Automatic Control}, 58(1):3--18, 2013.

\end{thebibliography}
\end{document}